\documentclass[preprint,12pt]{elsarticle}
\usepackage{graphics}
\usepackage{graphicx}
\usepackage{subfigure}

\usepackage{latexsym}
\usepackage{amsmath}
\usepackage{graphicx}
\usepackage{times}
\usepackage{graphicx,multicol,enumerate,subfigure,amsmath}

\usepackage{fullpage}
\usepackage{setspace}
\usepackage[english]{babel}
\usepackage[version=3]{mhchem}
\usepackage{amsmath, amssymb, amsthm}
\usepackage{verbatim, latexsym}
\usepackage{colortbl}

\newtheorem{theorem}{Theorem}
\newtheorem{remark}[theorem]{Remark}

\usepackage{soul}
\setstcolor{red}

\biboptions{sort&compress}
\journal{Journal of Computational Physics}

\onehalfspace

\begin{document}


\title{Multiscale Empirical Interpolation for Solving Nonlinear PDEs
  using Generalized Multiscale Finite Element Methods}

\author{\textbf{Victor M. Calo}$^{1,2}$, \textbf{Yalchin
    Efendiev}$^{1,3}$, \textbf{Juan Galvis}$^4$}

\author{\textbf{Mehdi Ghommem}$^1$ \corref{cor1}}
\cortext[cor1]{Email address : mehdig@vt.edu}

\address{$^{1}$
  Center for Numerical Porous Media (NumPor) \\
  King Abdullah University of Science and Technology (KAUST) \\
  Thuwal 23955-6900, Kingdom of Saudi Arabia}

\address{$^{2}$ Applied Mathematics \& Computational Science and Earth
  Sciences \& Engineering \\
  King Abdullah University of Science and Technology (KAUST) \\
  Thuwal 23955-6900, Kingdom of Saudi Arabia}

\address{$^{3}$ Department of Mathematics \& Institute for Scientific Computation (ISC) \\
  Texas A\&M University \\
  College Station, Texas, USA}

\address{$^{4}$
  Departamento de Matem\'aticas, Universidad Nacional de Colombia,\\
  Carrera 45 No 26-85 - Edificio Uriel Gutierr\'ez, Bogot\'a D.C. -
  Colombia}

\begin{abstract}

  In this paper, we propose a multiscale empirical interpolation
  method for solving nonlinear multiscale partial differential
  equations.  The proposed method combines empirical interpolation
  techniques and local multiscale methods, such as the Generalized
  Multiscale Finite Element Method (GMsFEM).  To solve nonlinear
  equations, the GMsFEM is used to represent the solution on a coarse
  grid with multiscale basis functions computed offline.  Computing
  the GMsFEM solution involves calculating the residuals on the fine
  grid.  We use empirical interpolation concepts to evaluate the
  residuals and the Jacobians of the multiscale system with a
  computational cost which is proportional to the coarse scale problem
  rather than the fully-resolved fine scale one.  Empirical
  interpolation methods use basis functions and an inexpensive
  inversion which are computed in the offline stage for finding the
  coefficients in the expansion based on a limited number of nonlinear
  function evaluations.  The proposed multiscale empirical
  interpolation techniques: (1) divide computing the nonlinear
  function into coarse regions; (2) evaluate contributions of
  nonlinear functions in each coarse region taking advantage of a
  reduced-order representation of the solution; and (3) introduce
  multiscale proper-orthogonal-decomposition techniques to find
  appropriate interpolation vectors.  We demonstrate the effectiveness
  of the proposed methods on several examples of nonlinear multiscale
  PDEs that are solved with Newton's methods and fully-implicit time
  marching schemes.  Our numerical results show that the proposed
  methods provide a robust framework for solving nonlinear multiscale
  PDEs on a coarse grid with bounded error.

\end{abstract}

\maketitle

\section{Introduction}

Solving nonlinear Partial Differential Equations (PDEs) with multiple
scales and/or high-contrast in media properties is computationally
expensive because of the disparity between scales that need to be
represented and the inherent nonlinearities. For this reason,
coarse-grid computational models are often used.  These include
Galerkin multiscale finite elements~\cite{Arbogast_two_scale_04,
  Chu_Hou_MathComp_10,ee03, egw10,eh09,ehg04, GhommemJCP2013}, mixed
multiscale finite element methods~\cite{aarnes04, ae07,
  Arbogast_Boyd_06,Iliev_MMS_11}, the multiscale finite volume
method~\cite{jennylt03}, mortar multiscale
methods~\cite{Arbogast_PWY_07, Wheeler_mortar_MS_12}, and variational
multiscale methods~\cite{hughes98}.  The coarse-grid models for
nonlinear PDEs can be divided into several classes.  One of them
includes constructing nonlinear operators that allow downscaling from
coarse-grid functions to fine-grid functions~\cite{ep04, ep05}. Other
types of approaches involve designing linear multiscale basis
functions (which are constructed using nonlinear PDEs) on the coarse
grid and using these basis functions for approximating the
solution~\cite{egh12}.

One of the challenges in these coarse-grid models is evaluating
nonlinear functionals for computing the residual and the Jacobian of
the operator.  Some of the widely used techniques are based on
Empirical Interpolation Methods (EIM)~\cite{barrault04, cs10, ohl12}.
In the offline stage, the response of the nonlinear function is
evaluated to yield solution snapshots and interpolation vectors are
constructed using the Proper Orthogonal Decomposition (POD).
Evaluating a few components of the nonlinear function allows
performing rapid evaluation of the nonlinear function in the online
stage with sufficient accuracy~\cite{cs10}.  This procedure allows for
rapid evaluations of the nonlinear function's response in the online
stage.  When solving multiscale equations, these empirical
interpolation techniques can be expensive because of the large problem
sizes. In this paper, we design an efficient multiscale empirical
interpolation framework.

Our proposed empirical interpolation techniques are based on
multiscale finite element approximations that are often used for
solving multiscale problems on coarse grids. The main idea of these
methods is to construct local multiscale basis functions for
approximating the solution over each coarse patch. Constructing these
basis functions uses offline coarse-grid spaces. These spaces are
constructed via judicious choices of snapshot spaces and local
spectral problems that are motivated by the analysis~\cite{egh12}.
The main objective of this paper is to design efficient local
multiscale empirical interpolations that can be used in conjunction
with Generalized Multiscale Finite Element Method (GMsFEM) to solve
nonlinear multiscale problems with a computational cost proportional
to the number of coarse-scale degrees of freedom rather than the fully
resolved mesh.

When GMsFEM is used for nonlinear problems, we evaluate the residual
and Jacobian in each Newton iteration which require fine-grid
calculations and incur a high computational cost.  We use POD-based
empirical interpolation and follow the Discrete Empirical
Interpolation Method (DEIM) introduced in~\cite{cs10}. For this
reason, we refer to our method as multiscale DEIM.  The main
ingredients of the proposed multiscale DEIM are:
\begin{itemize}
\item The evaluation of the nonlinear functions is performed on the
  subregions associated with those of the local multiscale basis
  functions and the global coupling of these local representations is
  performed.
\item In each subregion, an empirical interpolation is used as a
  representation of the nonlinear function which is computed based on
  a small dimensional multiscale representation of the solution space.
\item For multiscale high-contrast problems, interpolation vectors are
  computed using appropriate spectral problems which involve local
  multiscale basis functions and ensure convergence of the
  method. This modified spectral problems are constructed using bounds
  provided by the numerical analysis of the system.
\end{itemize}

In this paper, we investigate these issues and show that one can
design efficient empirical interpolation techniques in conjunction
with GMsFEM. We provide error estimates and the derivations of new
spectral problems for computing empirical interpolation vectors.

Representative numerical results are presented in the paper.  We
consider three examples: (1) nonlinear function evaluations using
GMsFEM basis functions; (2) steady-state multiscale elliptic equation
with nonlinear forcing; and (3) nonlinear multiscale parabolic
equations with nonlinear diffusion coefficients. In all examples, we
show that, by using a multiscale empirical interpolation, we can
obtain an accurate solution approximation at a cost that scales with
the coarse-grid size and is independent of the fine-grid size.  In
particular, a few interpolation modes are needed in all these
examples.

The paper is organized as follows. In Section~\ref{sec:review}, we
discuss empirical interpolation methods and local multiscale
techniques.  In Section~\ref{MsDEIM}, we introduce a multiscale
empirical interpolation and discuss its convergence. In Section
\ref{sec:app}, applications of the multiscale empirical interpolation
techniques to nonlinear PDEs are studied.  In this section, we discuss
GMsFEM and the use of multiscale interpolation techniques in Newton
methods. Numerical results are presented in Section~\ref{sec:num}.

\section{Review of basic concepts}
\label{sec:review}
\subsection{Discrete Empirical interpolation method}\label{sec:deim}

In this paper, we use the Discrete Empirical Interpolation Method
(DEIM)~\cite{cs10} for the local approximation of nonlinear functions
though other empirical interpolation methods can also be
used~\cite{barrault04}. DEIM approximates a nonlinear function by
means of an interpolatory projection of a few selected global
snapshots of the function. The idea is to represent a function over
the domain while using empirical snapshots and information in some
locations (or components).

We briefly review DEIM following~\cite{cs10}.  First, we give some
motivation for using DEIM. Let $ {f}(\tau) \in \mathbb{R}^{n}$ denote
a nonlinear function, where $\tau$ refers to any control parameter.
In a reduced-order modeling, the state vector $\tau$ is typically
assumed to have a reduced-order representation, i.e., $\tau\in
\mathbb{R}^{\widehat{n}}$ can be represented by fewer basis vectors,
$\zeta_1$,..., $\zeta_l$, where $l\ll \widehat{n}$.  Here, in general,
$\widehat{n}$ can be different from $n$, though we can assume
$\widehat{n}=n$. The reduced-order representation of $\tau$ usually
leads us to look for a reduced-order representation of the nonlinear
functions $f(\tau)$. The procedure of finding the reduced-order
representation of $f(\tau)$ consists of two steps.  In the first step,
we would like to find $m$ basis vectors (where $m$ is much smaller
than $n$), $\psi_1$,..., $\psi_m$, such that $f(\tau)$ can be
approximated in the space spanned by these vectors.  The error of this
approximation is given by POD error which represents approximation of
all possible $f(\tau_i)$'s in the space spanned by $\psi_1$,...,
$\psi_m$~\cite{cs10}.  In the second step, we identify a
reduced-dimensional linear system that allows finding the
representation of $f(\tau)$ in the space without involving $n$
operations. In this step, we typically identify $m$ equations that
allow finding the coordinates of $f(\tau)$ in the space spanned by
$\psi_1$,..., $\psi_m$.

We assume an approximation of the function
${f}$ obtained by projecting it into a subspace spanned by the basis
functions (snapshots) $ \Psi=( \psi_{_1},\cdots, \psi_{_m}) \in
\mathbb{R}^{n\times m}$ which are obtained by forward simulations.
We
write
\begin{equation}\label{fapp}
  {f}(\tau) \approx  \Psi  {d}(\tau).
\end{equation}
To compute the coefficient vector $ d$, we select $m$ rows of
(\ref{fapp}) and invert a reduced system to compute $d(\tau)$. This
can be formalized using the matrix $\mbox{P}$
$$
\mbox{P} = [ {e}_{\wp_1},\cdots, {e}_{\wp_m}] \in \mathbb{R}^{n\times
  m},
$$
where $ {e}_{\wp_i}=[0,\cdots,0,1,0,\cdots,0]^T\in \mathbb{R}^{n}$ is
the $\wp_i^{\mbox{th}}$ column of the identity matrix $\mbox{I}_n \in
\mathbb{R}^{n\times n}$ for $i=1,\cdots,m$. Multiplying
Equation~\eqref{fapp} by $ {P}^T $ and assuming that the matrix
$\mbox{P}^T \Psi$ is nonsingular, we obtain
\begin{equation}\label{fapp2}
  {f}(\tau) \approx  \Psi  d(\tau)= \Psi (\mbox{P}^T
  \Psi)^{-1}\mbox{P}^T {f}(\tau).
\end{equation}

To summarize, approximating the nonlinear function $ {f}(\tau)$, as
given by Equation~\eqref{fapp2}, requires the following:
\begin{itemize}
\item computing the projection basis $ \Psi=( \psi_{_1},\cdots,
  \psi_{_m}).$
\item identifying the indices $\{\wp_1,\cdots,\wp_m\}.$
\end{itemize}

To determine the projection basis $ \Psi=( \psi_{_1},\cdots,
\psi_{_m})$, we collect function evaluations in an $n \times n_s$
matrix $\mbox{F}=[ {f}(\tau_1),\cdots, {f}(\tau_{n_s})]$ and employ
the POD technique to select the most energetic modes.  This selection
uses the eigenvalue decomposition of the square matrix
$\mbox{F}^T\mbox{F}$ and selects the important modes using the
dominant eigenvalues.  These modes are used as the projection basis in
the approximation given by Equation~\eqref{fapp}.  In
Equation~\eqref{fapp2}, the term $ \Psi (\mbox{P}^T \Psi)^{-1} \in
\mathbb{R}^{n\times m}$ is computed once and stored.  The $d(\tau)$ is
computed using the values of the function $f(\tau)$ at $m$ points with
the indices $\wp_1,\cdots,\wp_m$ (identified using {the DEIM
  algorithm}). The computational saving is due to fewer evaluations of
$f(\tau)$. We refer to~\cite{cs10} for further details.

\begin{remark}[Inclusion of a priori multiscale
information on the POD selection procedure]\label{rem:m1m2}
In this paper, we explore ways to include the information about
heterogeneities in the POD selection process.  POD selects the
dominant eigenpairs of $(\mbox{F}^T\mbox{F}) z_i=\lambda_i z_i$, where
$z_i$ are coordinates in the space $\text{Span}\{ {f}(\tau_1),\cdots,
{f}(\tau_{n_s}) \}$.  In multiscale high-contrast problems, one needs
to take into account the heterogeneities when calculating the POD
modes.  Depending on the application in mind and on a priori
information on the space of functions $u$ for which $f(u)$ needs to be
computed, it is possible to improve the approximation using different
inner products, say, represented by $M_1$ and $M_2$, and perform a
spectral selection based on the eigenvalues of the modified eigenvalue
problem $\mbox{F}^TM_1\mbox{F}=\lambda \mbox{F}^TM_2\mbox{F}z$.
\end{remark}

\subsection{Local multiscale interpolation}

\begin{figure}
  \begin{center}
\includegraphics[width=0.5\textwidth]{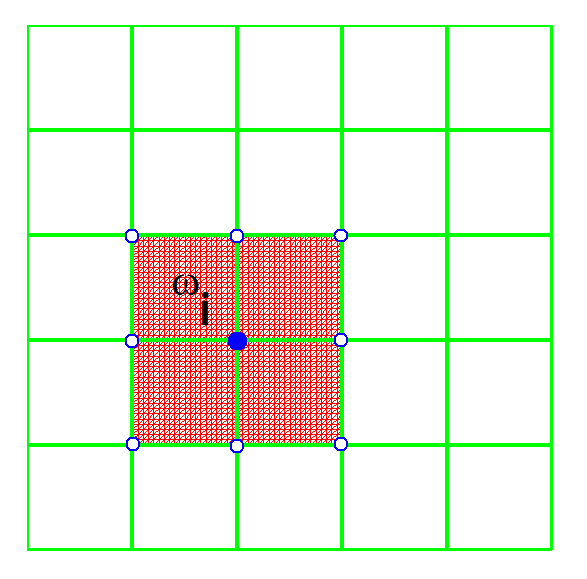}
  \end{center}
  \caption{A coarse-grid decomposition of a two-dimensional square
    domain.  We highlight the $i$-th coarse-node (blue circle), its
    neighborhood $\omega_i$ and all other coarse nodes in $\omega_i$
    (white circles). The fine-mesh is depicted only for coarse-grid
    blocks forming the highlighted support $\omega_i$ of the $i$-th
    coarse basis. }
  \label{fig:region}
\end{figure}

In multiscale problems, one needs to construct local approximations of
the solution using appropriately designed multiscale basis
functions. We use GMsFEM, a general local multiscale
strategy~\cite{egh12}.  These methods construct local multiscale basis
functions to approximate the solution over each coarse
patch. Constructing these basis functions uses local offline
spaces. These spaces are constructed via judicious
  choices of snapshot spaces
and local spectral problems that are motivated by the analysis.  Local
multiscale basis functions are defined on a coarse grid where each
coarse-grid block is a union of fine-grid blocks (see
Figure~\ref{fig:region} for a schematic representation of the coarse
and fine grids).

Next, we discuss the coarse-grid projection operators without going
into details. A detailed construction is presented in
Section~\ref{sec:gmsfem}.  Assume that the fine-scale problem has $n_f$
degrees of freedom and that the coarse-scale problem has $n_c$ degrees
of freedom. We design a matrix $\Phi$ of size $n_f\times n_c$ whose
transpose $\Phi^T$ of size $n_c\times n_f$ maps the fine-grid vectors to
vectors of coarse degrees of freedom (coarse vectors). We refer to
$\Phi$ as downscaling (or coarse-to-fine) operator and to $\Phi^T$ as
upscaling (or fine-to-coarse projection) operator. This construction
depends on several ingredients, such as the snapshot space and the
eigenvalues of the offline and online problems (see
Section~\ref{sec:gmsfem} for details).

\begin{figure}[htb]
  \begin{center}
\includegraphics[width=0.9\textwidth]{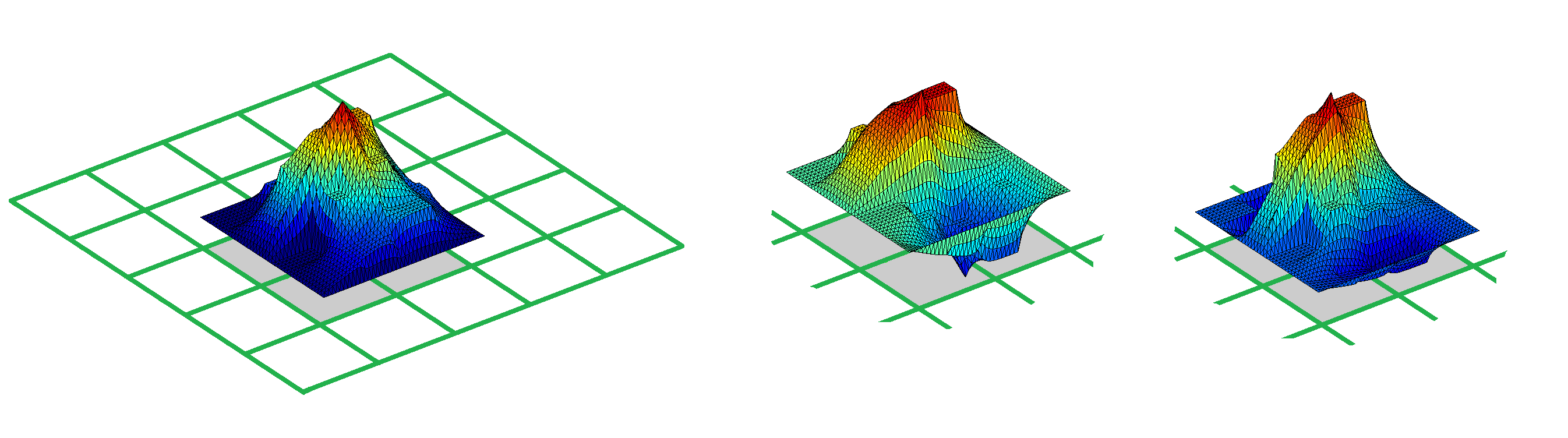}
  \end{center}
  \caption{Illustration of three basis functions supported in the
    highlighted neighborhood. }
  \label{fig:basis}
\end{figure}

Locality is key to the design of the new class of DEIM techniques
proposed herein.  Each column of $\Phi$ is a multiscale basis function
supported locally on a coarse region (see Figure~\ref{fig:basis} for
an example of a basis in a coarse region). One usually has several
multiscale basis functions per coarse-grid block.  In the design of
our multiscale DEIM, we identify multiscale basis functions (columns
of $\Phi$) which have support in $\omega_i$. Only these selected
multiscale basis functions are used in obtaining the local DEIM
approximation of $f(u)$.

When solving nonlinear PDEs, one writes the residual on the fine grid
as
\begin{equation}
  \label{eq:fine-res}
  R(u)=0,
\end{equation}
where $R(u)$ is the residual of nonlinear PDE and $u$ is the fine-grid
solution. For instance, in a nonlinear flow problem we may have a
nonlinear partial differential equation (in strong form) given by
$R(u)= \alpha(x) {\partial u \over \partial t} -\mbox{div} \big(
\kappa(x,u,\mu_\kappa) \, \nabla u \big)-g(u,x,\mu) $; see
Sections~\ref{Newton} and~\ref{sec:num} for further details and
examples.  Here, both $u$ and $R(u)$ are $n$-dimensional vectors
defined on a fine grid. Using the projection operator $\Phi$, we
project~\eqref{eq:fine-res} onto the coarse degrees of freedom (noting
that $\Phi z$ is an approximation of the fine-grid solution, that is,
$u\approx\Phi z$)
\begin{equation}
\label{eq:coarse-res}
\Phi^T R(\Phi z)=0.
\end{equation}
This equation is formulated on the coarse degrees of freedom
constructed on the coarse-grid; however, computing the residual $
R(\Phi z)$ requires fine-grid evaluations.  Moreover, computing the
Jacobians for each Newton iteration, defined as
\[
J(z)=\nabla_z R(\Phi z),
\]
also requires fine-grid evaluations. Here, our main goal is to use the
multiscale DEIM to compute $ R(\Phi z)$ and $J(z)$ efficiently. In
particular, using the multiscale DEIM approximation similar
to~\eqref{fapp}, we can write
\[
R(\Phi z)\approx \Psi d(z).
\]
Consequently, the residual computation involves
\begin{equation}
\label{eq:coarse-res1}
\Phi^T \Psi d(z),
\end{equation}
which can be efficiently computed by pre-computing $\Phi^T \Psi $. A
similar procedure can be applied to compute the Jacobian $J(z)$.

\section{Multiscale Discrete Empirical Interpolation
  Methods}\label{MsDEIM}

\subsection{Algorithm}

We detail the multiscale Discrete Empirical Interpolation Methods
(multiscale DEIM). A key idea is that, instead of solving the
fine-grid problem~\eqref{eq:fine-res}, we solve the coarse
problem~\eqref{eq:coarse-res}. An important issue in the solution of
the coarse problem~\eqref{eq:coarse-res} is the evaluation of the
nonlinear terms.  In order to fix ideas, in this paper we consider
scalar problems and the evaluation of the nonlinear terms of the from
$f(u)$ where (the fine-grid finite element function) $u:D\to
\mathbb{R}$ and $f:\mathbb{R}\to \mathbb{R}$.  In our approach we use
the DEIM procedure to efficiently evaluate the nonlinear terms.  We
stress the following main observations that are explored and
ultimately motivate the design of the multiscale DEIM procedure
presented below.

\begin{itemize}
\item In applications to multiscale PDEs (where one
  solves~\eqref{eq:coarse-res} instead of the fine-grid
  problem~\eqref{eq:fine-res}) the nonlinear functional $f$ needs to
  be evaluated with vectors of the form $u=\Phi z$ that are the
  downscaling of solutions obtained by reduced-order models. Thus,
  $f(\Phi z)$ needs to be computed in the span of coarse-grid snapshot
  vectors which has a reduced dimension.
\item Due to the fact that multiscale basis functions are supported on
  a coarse-grid neighborhood, it is sufficient to obtain the DEIM
  approximation in each coarse-node neighborhood.
\item More elaborate spectral selections may be needed to identify the
  elements of empirical interpolation vectors such that the resulting
  multiscale DEIM approximation is accurate in adequate norms that
  depend on physical parameters such as the contrast and small scales.
\end{itemize}

We recall that each column of $\Phi$ is the vector representation of a
(fine grid) finite element function with local support so we can
identify columns of $\Phi$ with locally supported basis functions.
Let us introduce the notation $I^{\omega_i} $ which represents the set
of indexes of coarse basis functions (which correspond to the columns
of $\Phi$) such that these basis functions have support on $\omega_i$
(see Figures~\ref{fig:region} and~\ref{fig:basis}).  Furthermore, we
introduce partition of unity diagonal matrices defined on $\omega_i$
denoted by $\mathcal{D}_i$ and such that
\begin{equation}
\label{eq:PU}
\sum_i \mathcal{D}_i = \mbox{I}_n,
\end{equation}
where $\mbox{I}_n$ is the identity matrix of size $n\times n$. Here,
for each node $i$, the $\mathcal{D}_i$ is a diagonal matrix with the
main diagonal that consists of a partition of unity vector $\chi_i$,
i.e.,
$\mathcal{D}_i=\chi_i(x_j)\delta_{ji}$, where $x_j$ are nodal
points and $\delta_{ji}$ is the Kronecker symbol.

We note that $f(u)$ can be written as
\[
f(u)=\sum_i \mathcal{D}_i f(u).
\]
Observe that $ \mathcal{D}_i f(u)$ is defined only on the fine-grid
nodes of $\omega_i$. Since $f:\mathbb{R}\to \mathbb{R}$, and thus only
a few basis functions (with indices from $I^{\omega_i}$ ) contribute
to it since for $u=Rz=\sum_{j=1}^{n_c} z_j\phi_j$ and $x\in \omega_i$
we have $u(x)=\sum_{j\in I^{\omega_i}} z_j \phi_j(x)$. Thus, we have
that
\[
f(\Phi z)=\sum_i \mathcal{D}_i f(\Phi z)=\sum_i \mathcal{D}_i
f(\sum_{j\in I^{\omega_i}} z_j \phi_j),
\]
where $\phi_j$ are basis vectors ($j$-th columns of $\Phi$).


Next, in each neighborhood $\omega_i$, we can perform an empirical
interpolation locally since we can write
\[
\sum_i \mathcal{D}_i f(\sum_{j\in I^{\omega_i}} z_j \phi_j) =
\sum_i \mathcal{D}_i f^{\omega_i}(\sum_{j\in I^{\omega_i}} z_j
\phi_j^{\omega_i})=
\sum_i \mathcal{D}_i f^{\omega_i}(\Phi^{\omega_i}z^{\omega_i}),
\]
where $f^{\omega_i}$ is the restriction of $f$ to $\omega_i$,
$\phi_j^{\omega_i}$ are the components of $\phi_j$ in $\omega_i$,
$\Phi^{\omega_i}$ is the matrix with columns $\phi_j^{\omega_i}$, and
$z^{\omega_i}$ represents a vector containing only entries $z_j$ with
$j\in I^{\omega_i}$. Therefore, noting that $f^{\omega_i}$ depends on
a few $z$'s, we perform an empirical interpolation using DEIM.  For
each coarse region $\omega_i$ we apply DEIM as introduced in
Section~\ref{sec:deim} to construct an approximation of the
type~\eqref{fapp} for the function $\tau \mapsto
f^{\omega_i}(\Phi^{\omega_i}\tau)$ with $\tau=z^{\omega_i}$.
According to~\eqref{fapp} we obtain $\Psi^{\omega_i}$ and
$d^{\omega_i}$ such that the following approximation holds,
\begin{equation}
\label{eq:deim_local}
f^{\omega_i}(\sum_{j\in I^{\omega_i}} z_j \phi_j^{\omega_i})\approx
\widetilde{f}^{\omega_i}(\sum_{j\in I^{\omega_i}} z_j
\phi_j^{\omega_i}):=\Psi^{\omega_i} d^{\omega_i}(z^{\omega_i}).
\end{equation}
With this empirical interpolation, we have
\begin{equation}
\label{eq:deim_local1}
f(\Phi z)\approx \widetilde{f}(\Phi z):=\sum_i \mathcal{D}_i
\widetilde{f}^{\omega_i}(\sum_{j\in I^{\omega_i}} z_j
\phi_j^{\omega_i})=\sum_i \mathcal{D}_i \Psi^{\omega_i}
d^{\omega_i}(z^{\omega_i}).
\end{equation}
In simulation, $\mathcal{D}_i \Psi^{\omega_i}$ can be
pre-computed and thus approximating $f(\Phi z)$ can be done at a lower
cost.

\begin{table}
\label{tab:DEIM_Algorithm}
  \begin{center}
    \begin{tabular}{l l}
      \hline
      \hline
      {DEIM Algorithm} & \emph{Multiscale Discrete Empirical Interpolation Method}   \\
      \hline
      {Input:} & Coarse grid, partition of unity matrices $\mathcal{D}_i$, multiscale basis matrix $\Phi^{\omega_i}$,\\
      & a projection basis matrix $ \Psi^{\omega_i}$
      obtained by applying POD on a sequence of\\
      & $n_s$ function evaluations \\
      \\
      & For each $\omega_i$\\
      &  {1.} Set $[\mid\rho\mid,\wp_1^{\omega_i}]=\max\{\mid \psi_{_1}^{\omega_i}\mid\}$ \\
      &  {2.} Set $ \Psi^{\omega_i}=( \psi_{_1}^{\omega_i})$, $\mbox{P}=( {e}_{\rho_1})$, , and $\overrightarrow{\wp}^{\omega_i}=(\wp_1^{\omega_i})$\\
      &  {3.}  {for} $\ell = 2,\cdots , m$  {do}  \\
      & $\;\;\;\;\;\;\;\;\;$ Solve $((\mbox{P}^{\omega_i})^T   \Psi^{\omega_i})  {d}^{\omega_i}=(\mbox{P}^{\omega_i})^T  \psi_{_{\ell }}^{\omega_i}$    \\
      & $\;\;\;\;\;\;\;\;\;$ Compute $ {r} =   \psi_{_{\ell }}^{\omega_i} -  \Psi^{\omega_i}  {d}^{\omega_i}$    \\
      & $\;\;\;\;\;\;\;\;\;$ Compute $[\mid\rho\mid,\wp_{\ell}^{\omega_i}]=\max\{\mid {r}\mid\}$   \\
      & $\;\;\;\;\;\;\;\;\;$ Set $ \Psi^{\omega_i}=( \Psi^{\omega_i}\;
      \psi_{_{\ell}}^{\omega_i})$,
      $\mbox{P}^{\omega_i}=(\mbox{P}^{\omega_i}\;
      {e}_{\rho_{\ell}^{\omega_i}})$, and
       $\wp^{\omega_i}=\left(
        \begin{array}{c}
          \wp^{\omega_i} \\
          \wp_{\ell}^{\omega_i} \\
        \end{array}
      \right)$
        \\

      & $\;\;\;\;$  {end for}       \\
                                                                                                                                               & $d^{\omega_i}=((P^{\omega_i})^T\Psi^{\omega_i})^{-1}\Psi^{\omega_i}$\\
                                                                                                                                               {Output:} & the interpolation indices $\overrightarrow{\wp}^{\omega_i}=(\wp_1^{\omega_i},\cdots,\wp_m^{\omega_i})^T$ and \\
                                                                                                                                               & $f(\Phi z)\approx \sum_i \mathcal{D}_i \Psi^{\omega_i} d^{\omega_i}(z^{\omega_i})$\\
                                                                                                                                               & \\
                                                                                                                                               & \\
                                                                                                                                               \hline
  \hline
\end{tabular}
\caption{Multiscale Discrete Empirical Interpolation
  Method.}
\end{center}
\end{table}

The algorithm for multiscale DEIM, used to compute the approximation
in~\eqref{eq:deim_local} for each coarse region, is presented in
Table~\ref{tab:DEIM_Algorithm}.

\subsection{Analysis}
In the rest of the section, we present the analysis of the multiscale
DEIM designed above.  In order to simplify the expressions, we use the
notation $A\preceq B$ to signify that there is a constant $C$ such
that $A\preceq C B$, where this constant $C$ is independent of the
vectors involved and of dimension $n$.

We start by using the definition of $\widetilde{f}(\Phi z)$
in~\eqref{eq:deim_local1} to get,
\begin{eqnarray}
  \nonumber
  \|f(\Phi z) - \widetilde{f}(\Phi z)\|^2&=&
  \|\sum_i \mathcal{D}_i f(\Phi z) - \widetilde{f}(\Phi z)\|^2\\
  \nonumber
  &=&
  \|\sum_i \mathcal{D}_i (f(\Phi z) -
  \widetilde{f}^{\omega_i}(\Phi^{\omega_i} z^{\omega_i}))\|^2 \\
  &\preceq&
  \nonumber
  \sum_i \|f(\Phi^{\omega_i} z^{\omega_i}) -
  \widetilde{f}^{\omega_i}(\Phi^{\omega_i} z^{\omega_i})\|^2 \\
  &=&\sum_i \|f^{\omega_i}(\Phi^{\omega_i} z^{\omega_i}) -
  \Psi^{\omega_i} d^{\omega_i}(z^{\omega_i})\|_{\omega_i}^2.
  \label{eq:analysis_step1}
\end{eqnarray}
Here, we have applied a triangle inequality and the fact that the
entries of the diagonal matrices $\mathcal{D}_i$ are less than one.
The hidden constant is the maximum number of coarse nodes in
$I^{\omega_i}$, that is, $\max_{i} \#(I^{\omega_i})$, where $\#$ is
the cardinality of the set.  This number is finite, is independent of
the fine-grid parameter $n$ and depends only on the coarse
triangulation configuration.

Following DEIM estimates in~\cite{cs10} (see Lemma 3.2.),
we can write
\begin{equation}
  \|f^{\omega_i}(\Phi^{\omega_i} z^{\omega_i}) -  \Psi^{\omega_i}
  d^{\omega_i}(z^{\omega_i})\|_{\omega_i}^2 \leq
  \|((P^{\omega_i})^T\Psi^{\omega_i})^{-1}\|^2  \|f^{\omega_i}  -
  \mathcal{P}_{\Psi^{\omega_i}}f^{\omega_i}\|_{\omega_i}^2.
\end{equation}
Here $((P^{\omega_i})^T\Psi^{\omega_i})^{-1}$ is assumed invertible
and $ \|f^{\omega_i} -
\mathcal{P}_{\Psi^{\omega_i}}f^{\omega_i}\|_{\omega_i}^2$ is the
projection of the error onto the space spanned by $\Psi^{\omega_i}$
(where $\mathcal{P}_{\Psi^{\omega_i}}$ is the orthogonal projection in
the vector norm $\| \cdot \|_{\omega_i}$ to the space spanned by
$\Psi^{\omega_i}$ ).

To estimate $ \|f^{\omega_i} -
\mathcal{P}_{\Psi^{\omega_i}}f^{\omega_i}\|_{\omega_i}^2$, we can
consider an error estimate that is typical for POD estimates. The
error over all snapshots is given by
\[
\sum_{j=1}^ {n^{\omega_i}}\|f^{\omega_i,j} -
\mathcal{P}_{\Psi^{\omega_i}}f^{\omega_i,j}\|_{\omega_i}^2
=\sum_{j=L_i+1}^{n^{\omega_i}}  \lambda^{\omega_i}_{j},
\]
where ${n^{\omega_i}}$ are the number of snapshots
in the region $\omega_i$: $\{f^{\omega_i,j}\}_{j=1}^{n^{\omega_i}}$. The $
\lambda^{\omega_i}_{j}$ are eigenvalues (in the decreasing order) of
$\Psi^{\omega_i}(\Psi^{\omega_i})^T$, and $L_i$ is the number of
selected eigenvalues whose corresponding eigenvectors are used in the
empirical interpolation.  This estimate requires that the snapshots
used in the above estimate appear in the local POD. One can also show
that, for an arbitrary element $f^{\omega_i}$ in the span of
snapshots, we have
\[
\|f^{\omega_i} -
\mathcal{P}_{\Psi^{\omega_i}}f^{\omega_i}\|_{\omega_i}^2\preceq
\lambda^{\omega_i}_{L_i+1}
 \|f^{\omega_i}\|_{\omega_i}^2.
\]
With this latter estimate, we can write
\begin{eqnarray}
\nonumber
 \sum_i \|f^{\omega_i}(\Phi^{\omega_i} z^{\omega_i}) -
 \Psi^{\omega_i} d^{\omega_i}(z^{\omega_i})\|_{\omega_i}^2&\preceq&
\sum_i \|((P^{\omega_i})^T\Psi^{\omega_i})^{-1}\|^2 \|f^{\omega_i} -
\mathcal{P}_{\Psi^{\omega_i}}f^{\omega_i}\|_{\omega_i}^2\\
\nonumber
&\preceq&\sum_i\|((P^{\omega_i})^T\Psi^{\omega_i})^{-1}\|^2
\|f^{\omega_i} -
\mathcal{P}_{\Psi^{\omega_i}}f^{\omega_i}\|_{\omega_i}^2\\
\label{eq:analysis_step2}
&\preceq& \sum_i \|((P^{\omega_i})^T\Psi^{\omega_i})^{-1}\|^2
\lambda^{\omega_i}_{L_i+1}
 \|f^{\omega_i}\|_{\omega_i}^2.
\end{eqnarray}
By inserting~\eqref{eq:analysis_step2} into~\eqref{eq:analysis_step1},
we obtain the estimate
\begin{equation}
\|f(\Phi z) - \widetilde{f}(\Phi z)\|^2=
\max_i\{\|((P^{\omega_i})^T\Psi^{\omega_i})^{-1}\|^2
\lambda^{\omega_i}_{L_i+1}\} \|f\|^2.
\end{equation}
One may need a special inner product for multiscale high-contrast
equations to ensure an appropriate bound, as discussed in
Section~\ref{sec:MsPOD}.

\section{Applications to multiscale PDEs}
\label{sec:app}

In this section, we describe the offline-online computational
procedure that is used to solve the forward problem on a coarse grid.
We elaborate on possible choices for the associated bilinear
forms to be used in the coarse space construction.  Below we offer a
general outline of the multiscale procedure.  For details on the
constructions and further considerations we refer to~\cite{egh12,
  egt11, eglp13oversampling, eglmsMSDG}.
\begin{itemize}
\item[1.]  Offline computations:
\begin{itemize}
\item 1.0. Coarse grid generation.
\item 1.1. Construction of snapshot space used to compute the offline
  space.
\item 1.2. Construction of a small dimensional offline space by
  dimensional reduction in the space of local snapshots.
\end{itemize}
\item[2.] Online computations:
\begin{itemize}
\item 2.1. For each input parameter set, compute multiscale basis
  functions.
\item 2.2. Solution of a coarse-grid problem for given forcing term and
  boundary conditions.
\end{itemize}
\end{itemize}

\subsection{Problem setup}

We consider various non-linear elliptic equations of the form
\begin{equation}
\label{eq:original} \alpha(x) {\partial u \over \partial
    t} -\mbox{div} \big( \kappa(x,u,\mu_\kappa) \, \nabla u
  \big)=g(u,x,\mu) \, \, \text{in} \, D,
\end{equation}
where $u=0$ on $\partial D$. Performing multiscale simulations
requires appropriate local multiscale basis functions.  We discuss
this procedure next.  The procedure below identifies local basis
functions. Denoting these basis functions by $\phi_j$, we seek the
solution
\begin{equation}\label{umultiscale}
  u(x,t_n)=\sum_{i} z_i^{n} \phi_i
\end{equation}
that solves~\eqref{eq:original}. We employ an implicit time
discretization and use the Newton method to solve the resulting
nonlinear system at each time level. In this case, the residual form
of problem~\eqref{eq:original} can be written as
\begin{align}
\label{eq:residual12}
R_m(z^{n+1})&=\sum_j (z_j^{n+1}-z_j^n)\int_D \alpha(x) \phi_j
\phi_m\nonumber\\
&+\sum z_j^{n+1}\int_D \kappa(x,\sum_l
z_l^{n+1}\phi_l,\mu_\kappa)\nabla \phi_j \nabla\phi_m\nonumber\\&-
\int_D g(\sum_l z_l^{n+1}\phi_l,x,\mu)=0,
\end{align}
where $z^n$ is the solution at the previous time step and $z^{n+1}$ is
the value of the solution at the latest iteration level.

\subsection{Multiscale spatial discretization via GMsFEM}
\label{sec:gmsfem}

In the offline computation, we first construct a snapshot space
$V_{\text{snap}}^{\omega_i}$.  Constructing the snapshot space may
involve solving various local problems for different choices of input
parameters or different fine-grid representations of the solution in
each coarse region.  We denote each snapshot vector (listing the
solution at each node in the domain) using a single index and create
the following matrix
$$
\Phi_{\text{snap}} = \left[ \phi_{1}^{\text{snap}}, \ldots,
  \phi_{M_{\text{snap}}}^{\text{snap}} \right],
$$
where $\phi_j^{\text{snap}}$ denotes the snapshots and
$M_{\text{snap}}$ denotes the total number of functions to keep in the
local snapshot matrix construction.

In order to construct an offline space $V_{\text{off}}$, we perform a
dimension reduction process in the space of snapshots using an
auxiliary spectral decomposition. The main objective is to use the
offline space to efficiently (and accurately) construct a set of
multiscale basis functions to be used in the online stage. More
precisely, we seek for a subspace of the snapshot space such that it
can approximate any element of the snapshot space in the appropriate
sense which is defined via auxiliary bilinear forms. At the offline
stage the bilinear forms are chosen to be
\emph{parameter-independent}, such that there is no need to
reconstruct the offline space for each $\nu$ value, where $\nu$ is
assumed to be a parameter that represents $u$ and $\mu_\kappa$ in
$\kappa(x,u,\mu_\kappa)$.  For constructing the offline space, we use
the average of over the
coarse region $\omega_i$ in $\kappa(x,\nu)$.  Thus, $\nu$ represents
both the average of $u$ and $\mu$.  We consider the following
eigenvalue problem in the space of snapshots:
\begin{eqnarray}
  A^{\text{off}} \Phi_k^{\text{off}} &=& \lambda_k^{\text{off}}
  S^{\text{off}}\Phi_k^{\text{off}},  \label{offeig1}
\end{eqnarray}
where
\begin{align}
  \displaystyle A^{\text{off}}&= [a^{\text{off}}_{mn}] =
  \int_{\omega_i} \overline{\kappa}(x, \nu) \nabla
  \phi_m^{\text{snap}}
  \cdot \nabla \phi_n^{\text{snap}} = \Phi_{\text{snap}}^T \overline{A} \Phi_{\text{snap}},\label{eq:Abar}\\
  \displaystyle S^{\text{off}} &= [s^{\text{off}}_{mn}] =
  \int_{\omega_i} \widetilde{\overline{\kappa}}(x, \nu)
  \phi_m^{\text{snap}} \phi_n^{\text{snap}} = \Phi_{\text{snap}}^T
  \overline{S} \Phi_{\text{snap}}.\label{eq:Sbar}
\end{align}
In the definitions of $A^{\text{off}}$ and $S^{\text{off}} $ above,
the coefficient $\overline{\kappa}(x, \nu) $ is defined as a
parameter-averaged coefficient. The coefficient
$\widetilde{\overline{\kappa}}(x, \nu)$ can be chosen simply as
$\widetilde{\overline{\kappa}}(x, \nu)=\overline{\kappa}(x, \nu)$ or
in more sophisticated manners that include information about
multiscale finite element basis functions; we refer to~\cite{egh12}
for details and examples. In Equation~\eqref{eq:Abar}, $\overline{A}$
(similarly for $\overline{S}$ in~\eqref{eq:Sbar}) denotes a fine-scale
matrix, except that parameter-averaged coefficients are used in its
construction, and also that $\overline{A}$ is constructed by
integrating only on $\omega_i$.  To generate the offline space, we
then choose the smallest $M_{\text{off}}$ eigenvalues from
Equation~\eqref{offeig1} and form the corresponding eigenvectors in
the respective space of snapshots by setting $\phi_k^{\text{off}} =
\sum_j \Phi_{kj}^{\text{off}} \phi_j^{\text{snap}}$ (for $k=1,\ldots,
M_{\text{off}}$), where $\Phi_{kj}^{\text{off}}$ are the coordinates
of the vector $\Phi_{k}^{\text{off}}$. See~\cite{egh12,egw10} for
further details. We then create the offline matrices
 $$
 \Phi_{\text{off}} = \left[ \phi_{1}^{\text{off}}, \ldots,
   \phi_{M_{\text{off}}}^{\text{off}} \right]
$$
to be used in the online space construction.

The online coarse space is used within the finite element framework to
solve the original global problem, where continuous Galerkin
multiscale basis functions are used to compute the global solution. In
particular, we seek a subspace of the respective offline space such
that it can approximate well any element of the offline space in an
appropriate metric. At the online stage, the bilinear forms are chosen
to be \emph{parameter-dependent}.  The following eigenvalue problems
are posed in the offline space:
\begin{eqnarray}
  A^{\text{on}}(\nu) \Phi_k^{\text{on}} &=& \lambda_k^{\text{on}}
  S^{\text{on}}(\nu)\Phi_k^{\text{on}},  \label{oneig1}
\end{eqnarray}
where
\begin{align*}
  \displaystyle A^{\text{on}}(\nu) &= [a^{\text{on}}(\nu)_{mn}] =
  \int_{\omega_i} \kappa(x, \nu) \nabla \phi_m^{\text{off}} \cdot
  \nabla \phi_n^{\text{off}} = \Phi_{\text{off}}^T A(\nu)
  \Phi_{\text{off}},\\
\displaystyle S^{\text{on}} &= [s^{\text{on}}_{mn}] =
  \int_{\omega_i} \widetilde{{\kappa}}(x, \nu)
  \phi_m^{\text{off}} \phi_n^{\text{off}} = \Phi_{\text{off}}^T
  {S}(\nu) \Phi_{\text{off}},
 \end{align*}
 and $\kappa(x, \nu)$ and $\widetilde{{\kappa}}(x, \nu)$ are now
 parameter dependent.  As before, the coefficient
 $\widetilde{{\kappa}}(x, \nu)$ can be chosen simply as
 $\widetilde{{\kappa}}(x, \nu)={\kappa}(x, \nu)$ or can include the
 multiscale finite element basis functions.  The choice of
 $\widetilde{\kappa}$ has implications on the dimensions of the
 resulting coarse spaces; see~\cite{egh12}.  To generate the online
 space, we then choose the smallest $M_{\text{on}}$ eigenvalues
 from~\eqref{oneig1} and form the corresponding eigenvectors in the
 offline space by setting $\phi_k^{\text{on}} = \sum_j
 \Phi_{kj}^{\text{on}} \phi_j^{\text{off}}$ (for $k=1,\ldots,
 M_{\text{on}}$), where $\Phi_{kj}^{\text{on}}$ are the coordinates of
 the vector $\Phi_{k}^{\text{on}}$.  If $\kappa(x,u)=k_0(x) b(u)$,
 then one can use the parameter-independent case of GMsFEM. In this
 case, there is no need to construct the online space (or the online
 space is the same as the offline space). From now on, we denote the
 online space basis functions by $\phi_i$.

\subsection{Newton method and Newton-DEIM}\label{Newton}
We consider a time-dependent nonlinear flow governed by the following
parabolic partial differential equation
\begin{eqnarray}\label{parabolic}
  \frac{\partial u}{\partial t} -\mbox{div} \big( \kappa(x,u,\mu) \,
  \nabla u  \big) = h(x)   \, \, \text{in} \, D.
\end{eqnarray}
The finite element discretization of Equation~\eqref{parabolic} yields
a system of ordinary differential equations given by
\begin{eqnarray}
  \mbox{M} \dot{\mbox{U}} + \mbox{F}(\mbox{U}) = \mbox{H},
\end{eqnarray}
where
$$
\mbox{U} = \left(
               \begin{array}{cccc}
                 u_1 & u_2 & \cdots & u_{n_f} \\
               \end{array}
             \right)^T
$$
is the vector collecting the pressure values at the fine-scale nodes
with $n_f$ being the total number of fine-scale nodes and $
\mbox{H}$ is the right-hand-side vector obtained by discretization.
In our derivations and simulations, we assume that
\begin{equation}
\label{eq:defQ}
\kappa(x,u,\mu)= \sum_{q=1}^Q \kappa_q(x) b_q(u,\mu)
\end{equation}
for some $Q$. In general, $\kappa(x,u,\mu)$, can be approximated as
in~\eqref{eq:defQ} using the offline basis functions
(see~\cite{barrault04, egh12} for dealing with the nonlinearity in $u$
within each coarse region).  This results in
\begin{equation*}
\mbox{F}(\mbox{U},\mu) = \sum_{q=1}^Q\mbox{A}_q
\Lambda_1^q(\mbox{U},\mu) \mbox{U},
\end{equation*}
 where we have
\begin{align*}
\mbox{A}_{ij}^q&=\int_D \kappa_q \nabla \phi^0_i \cdot \nabla \phi^0_j
, \;\;\; \mbox{M}_{ij}=\int_D  \phi^0_i \phi^0_j\;,
\;\;\;\mbox{H}_{i}=\int_D  \phi^0_i h \;, \\
\Lambda_1^q(\mbox{U},\mu) &= \mbox{diag}\left(
               \begin{array}{cccc}
                 b_q(u_1,\mu) & b_q(u_2,\mu) & \cdots &
                 b_q(u_{N_f},\mu) \\
               \end{array}
             \right),
\end{align*}
and $\phi^0_i$ are piecewise linear basis functions defined on a fine
triangulation of $D$.

Employing the backward Euler scheme for the time
marching process, we obtain
\begin{eqnarray}
  \mbox{U}^{n+1}+\Delta t \;\mbox{M}^{-1}\mbox{F}(\mbox{U}^{n+1})=
  \mbox{U}^{n} + \Delta t \;\mbox{M}^{-1} \mbox{H},
\end{eqnarray}
where $\Delta t$ is the time-step size and the superscript $n$ refers
to the temporal level of the solution. We let
\begin{eqnarray}
  \mbox{R}(\mbox{U}^{n+1}) = \mbox{U}^{n+1}- \mbox{U}^{n} + \Delta t
  \;\mbox{M}^{-1}\mbox{F}(\mbox{U}^{n+1}) - \Delta t \;\mbox{M}^{-1}
  \mbox{H}
\end{eqnarray}
with derivative
\begin{align}
\mbox{J}(\mbox{U}^{n+1}) = D\mbox{R}(\mbox{U}^{n+1}) &= I + \Delta t
\;\mbox{M}^{-1} D \mbox{F}(\mbox{U}^{n+1}) \nonumber \\
&= I + \sum_{q=1}^Q \Delta t \;\mbox{M}^{-1}
\mbox{A}_q\Lambda_1^q(\mbox{U}^{n+1}) +\sum_{q=1}^Q \Delta t
\;\mbox{M}^{-1} \mbox{A}_q\Lambda_2^q(\mbox{U}^{n+1}),
\end{align}
where
\begin{equation*}
\Lambda_2^q(\mbox{U},\mu) = \mbox{diag}\left(
               \begin{array}{cccc}
                 \frac{\partial b_q(u_1,\mu)}{\partial u} &
                 \frac{\partial b_q(u_2,\mu)}{\partial u} &
                   \cdots &
                  \frac{\partial b_q(u_{N_f},\mu)}{\partial u}
               \end{array}
             \right),
\end{equation*}
and $D$ is the multi-variate gradient operator defined as
$D\mbox{R}(\mbox{U})=\partial \mbox{R}_i/\partial \mbox{U}_j$.
The scheme involves, at each time step, the following iterations
\begin{align}
  \mbox{J}(\mbox{U}^{n+1}_{(k)}) \Delta \mbox{U}^{n+1}_{(k)}&=
  -\Big(\mbox{U}^{n+1}_{(k)}- \mbox{U}^{n} + \Delta t
  \;\mbox{M}^{-1}\mbox{F}(\mbox{U}^{n+1}_{(k)}) - \Delta t
  \;\mbox{M}^{-1} \mbox{H}\Big) \nonumber \\
  \mbox{U}^{n+1}_{(k+1)}&=\mbox{U}^{n+1}_{(k)} + \Delta
  \mbox{U}^{n+1}_{(k)}, \nonumber
\end{align}
where the initial guess is $\mbox{U}^{n+1}_{(0)}=\mbox{U}^{n}\nonumber
$ and $k$ is the iteration counter. The above iterations are
repeatedly applied until $\parallel \Delta
\mbox{U}^{n+1}_{(k)} \parallel$ is less than a specific tolerance.

In our simulations, we use $Q=1$ (see~\eqref{eq:defQ}) for the
definition of $Q$) as our focus is on localized multiscale interpolation
of nonlinear functionals that arise in discretization of multiscale
PDEs.  With this choice, we do not need to compute the online
multiscale space (i.e., the online space is the same as the offline
space).

We use the solution expansion given by Equation~\eqref{umultiscale}
and employ the multiscale framework to obtain a set of $n_c$ ordinary
differential equations that constitute a reduced-order model; that is,
\begin{equation}\label{rom}
  \dot{z} = - (\Phi^T\mbox{M}\Phi)^{-1} \Phi^T \mbox{F} (\Phi z) +
  (\Phi^T\mbox{M}\Phi)^{-1} \Phi^T \mbox{H}.
\end{equation}
Thus, the original problem with $n_f$ degrees of freedom is reduced to a
dynamical system with $n_c$ dimensions where $n_c \ll n_f$.

The nonlinear term $(\Phi^T\mbox{M}\Phi)^{-1} \Phi^T \mbox{F} (\Phi
z)$ in the reduced-order model, given by Equation~\eqref{rom}, has a
computational complexity that depends on the dimension of the full
system $n_f$. This nonlinear term requires matrix multiplications and
full evaluation of the nonlinear function $\mbox{F}$ at the
$n_f$-dimensional vector $\Phi z(t)$. As such, solving the reduced
system still requires extensive computational resources and time. To
reduce this computational requirement, we use multiscale DEIM as
described in the previous section. In this case, computational savings
can be obtained in a forward run of the nonlinear model.

To solve the reduced system, we employ the backward Euler
scheme; that is,
\begin{eqnarray}
z^{n+1}+\Delta t \;\widetilde{\mbox{M}}^{-1}
\widetilde{\mbox{F}}(z^{n+1})= z^{n} + \Delta t
\;\widetilde{\mbox{M}}^{-1} \widetilde{\mbox{H}},
\end{eqnarray}
where $\widetilde{\mbox{M}}= \Phi^T\mbox{M}\Phi$,
$\widetilde{\mbox{F}}(z)=\Phi^T\mbox{F}(\Phi z)$, and
$\widetilde{\mbox{H}}=\Phi^T \mbox{H}$. We let
\begin{eqnarray}
\label{eq:msres}
\widetilde{\mbox{R}}(z^{n+1}) = z^{n+1}- z^{n} +\Delta t
\;\widetilde{\mbox{M}}^{-1} \widetilde{\mbox{F}}(z^{n+1}) - \Delta t
\;\widetilde{\mbox{M}}^{-1} \widetilde{\mbox{H}}
\end{eqnarray}
with derivative
\begin{align}
\widetilde{\mbox{J}}(z^{n+1})=D\widetilde{\mbox{R}}(z^{n+1}) &= I +
\Delta t \;\widetilde{\mbox{M}}^{-1} D \widetilde{\mbox{F}}(z^{n+1})
\nonumber \\
&= I + \sum_{q=1}^Q\Delta t \;\widetilde{\mbox{M}}^{-1}
\Phi^T\mbox{A}_q\Lambda_1^q(\Phi z^{n+1})\Phi + \sum_{q=1}^Q \Delta t
\;\widetilde{\mbox{M}}^{-1} \Phi^T\mbox{A}_q\Lambda_2^q(\Phi
z^{n+1})\Phi. \nonumber
\end{align}
The scheme involves, at each time step, the following iterations
\begin{align}\label{scheme}
  \widetilde{\mbox{J}}(z^{n+1}_{(k)}) \Delta z^{n+1}_{(k)}&=
  -\Big(z^{n+1}_{(k)}- z^{n} + \Delta t
  \;\widetilde{\mbox{M}}^{-1}\widetilde{\mbox{F}}(z^{n+1}_{(k)}) -
  \Delta t \;\widetilde{\mbox{M}}^{-1} \widetilde{\mbox{H}} \Big) \\
  z^{n+1}_{(k+1)}&=z^{n+1}_{(k)} + \Delta z^{n+1}_{(k)},
\end{align}
where the initial guess is $z^{n+1}_{(0)}=z^{n}\nonumber $. The above
iterations are repeated until $\parallel \Delta
z^{n+1}_{(k)} \parallel$ is less than a specific tolerance. We use
multiscale DEIM as detailed in Section~\ref{MsDEIM} to approximate the
nonlinear functions that appear in the residual $\widetilde{\mbox{R}}$
and the Jacobian $\widetilde{\mbox{J}}$ and, therefore, reduce the
number of function evaluations.

\subsection{The use of multiscale POD for DEIM}
\label{sec:MsPOD}

When designing empirical interpolation methods for high-contrast
problems, it is important to take into account heterogeneities and
perform an interpolation using special weighted norms. These norms are
derived based on the analysis discussed next.

Consider the approximation of
\begin{equation}
\label{eq:res1}
\sum_i c_i \int_D \kappa(x) b(\sum_l z_l \phi_l) \nabla \phi_i \nabla
\phi_j
\end{equation}
that appears in the residual function~\eqref{eq:msres}.  For
multiscale DEIM, we split this integral using partition of unity
matrices as discussed earlier in~\eqref{eq:PU} and re-write it for
the discretization of the stiffness matrix. We introduce the partition
of unity $\{\chi^{\omega_k}\}$ subordinated to the coarse regions
$\{\omega_k\}$.  The relation with the partition of unity matrices
$\mathcal{D}_k$ is that the action of $\mathcal{D}_k$ is the
multiplication of the corresponding finite element function by the
partition of unity function $\chi^{\omega_k}$.\label{def:xi} We obtain,

\begin{align}
\label{eq:EIM-PDE}
\sum_i c_i \int_D \kappa(x) b(\sum_l z_l \phi_l) \nabla \phi_i \nabla
\phi_j&=
\sum_{i,k} c_i \int_{\omega_k} \chi^{\omega_k}(x)\kappa(x)
b\left(\sum_l z_l \phi_l^{\omega_k}\right) \nabla \phi_i^{\omega_k}
\nabla \phi_j^{\omega_k}\nonumber\\
&=\sum_{i,k} c_i \int_{\omega_k} \chi^{\omega_k}(x) \kappa(x)
b\left(\sum_{l\in I^{\omega_k}} z_l \phi_l^{\omega_k}\right) \nabla
\phi_i^{\omega_k} \nabla \phi_j^{\omega_k}.
\end{align}
The interpolation for $b(\sum_{l\in I^{\omega_k}} z_l \phi_l^{\omega_k})$
is performed in the region $\omega_k$. Each term in the sum
(\ref{eq:EIM-PDE})
can be written as
\begin{equation}
\label{eq:EIM-weight}
\int_{\omega_k} m_{ij}^{\omega_k}(x) b(\sum_{l\in I^{\omega_k}} z_l
\phi_l^{\omega_k}),
\end{equation}
where $m_{ij}^{\omega_k}(x) = \chi^{\omega_k}(x) \kappa(x) \nabla
\phi_i^{\omega_k} \nabla \phi_j^{\omega_k}$.

To evaluate~\eqref{eq:EIM-weight}, one needs to take into account the
weight $m_{ij}^{\omega_k}(x) $ for high-contrast multiscale problems,
where $m_{ij}^{\omega_k}(x) $ can have very high values in some
subregions within $\omega_k$.  Otherwise, the accuracy of the method
substantially deteriorates.  However, the weighting function is not
uniquely defined and depends on particular basis functions that are
involved in the integration. We propose the use of a weight function
that is computed by summing over all indices $i,j$ (which provides an
upper bound for every $m_{ij}^{\omega_k}(x) $) and using only a few
dominant modes.  For high-contrast problems, $m_{ij}^{\omega_k}(x) $
is very high in the regions of high conductivity and this high value
can be estimated using only the first basis function.  We propose
using
\[
\widetilde{m}^{\omega_k}=\sum_{i} \kappa(x)|\nabla \phi_1^{\omega_i}|^2
\]
which is defined in the entire domain.  Here, $\phi_1^{\omega_k}$ is
the first basis function (associated to the first dominant mode in the
online GMsFEM selection procedure).  One can also use additional basis
functions to construct the weight $\widetilde{m}^{\omega_k}$.  For
implementing this procedure, we need to define $\Psi^{\omega_i}$ in
(\ref{eq:deim_local}) using weighted POD modes.  This can simply be
done using the matrix corresponding to the mass matrix
$\widetilde{m}^{\omega_k}$ as $M_1$ in Remark~\ref{rem:m1m2} ($M_2$ is
the identity matrix of appropriate size).

Using the above argument, one can write the error in the residuals of
multiscale DEIM $\widetilde{R}^D$ and $R(z)$.
\begin{equation}
\begin{split}
  R(z)-\widetilde{R}^D(z)= \sum_{i,k} c_i \int_{\omega_k}
  \chi^{\omega_k}(x) \kappa(x) (b(\sum_{l\in I^{\omega_k}} z_l
  \phi_l^{\omega_k}) - \sum_l d_l(z) \psi_l^{\omega_k} ) \nabla
  \phi_i^{\omega_k} \nabla \phi_j^{\omega_k}.
\end{split}
\end{equation}
Using an appropriate weight function in POD, we can control the error
$R(z)-\widetilde{R}^D(z)$ independently of the physical parameters
such as the contrast and small scales.  This error can be related,
under some conditions, to the error between the solution obtained
using multiscale DEIM and the solution obtained without DEIM.  If a
standard (instead of multiscale) POD is used for the selection of
empirical modes, then the error depends on the contrast and, in
particular, the error is proportional to the contrast. Similarly, the
error in the Jacobians corresponding to multiscale DEIM and without
DEIM can be controlled independently of the contrast if we choose the
weighted POD modes properly. Large errors in the Jacobians may lead
to a poor convergence of the Newton method as its convergence rate is
related to the quality of approximating the inverse of the
Jacobians. We have observed this in our numerical simulations.

The weighted POD procedure allows minimizing the residual error
between the snapshots of $b$'s (nonlinear function) and their
projections in the weighted norm. Consequently, the residual is mainly
minimized in the high-conductivity regions which can constitute a
small portion of the coarse block.

\section{Representative numerical experiments}
\label{sec:num}

\begin{figure}[htb]
  \begin{center}
    \subfigure[Case
    I]{\includegraphics[width=0.48\textwidth]{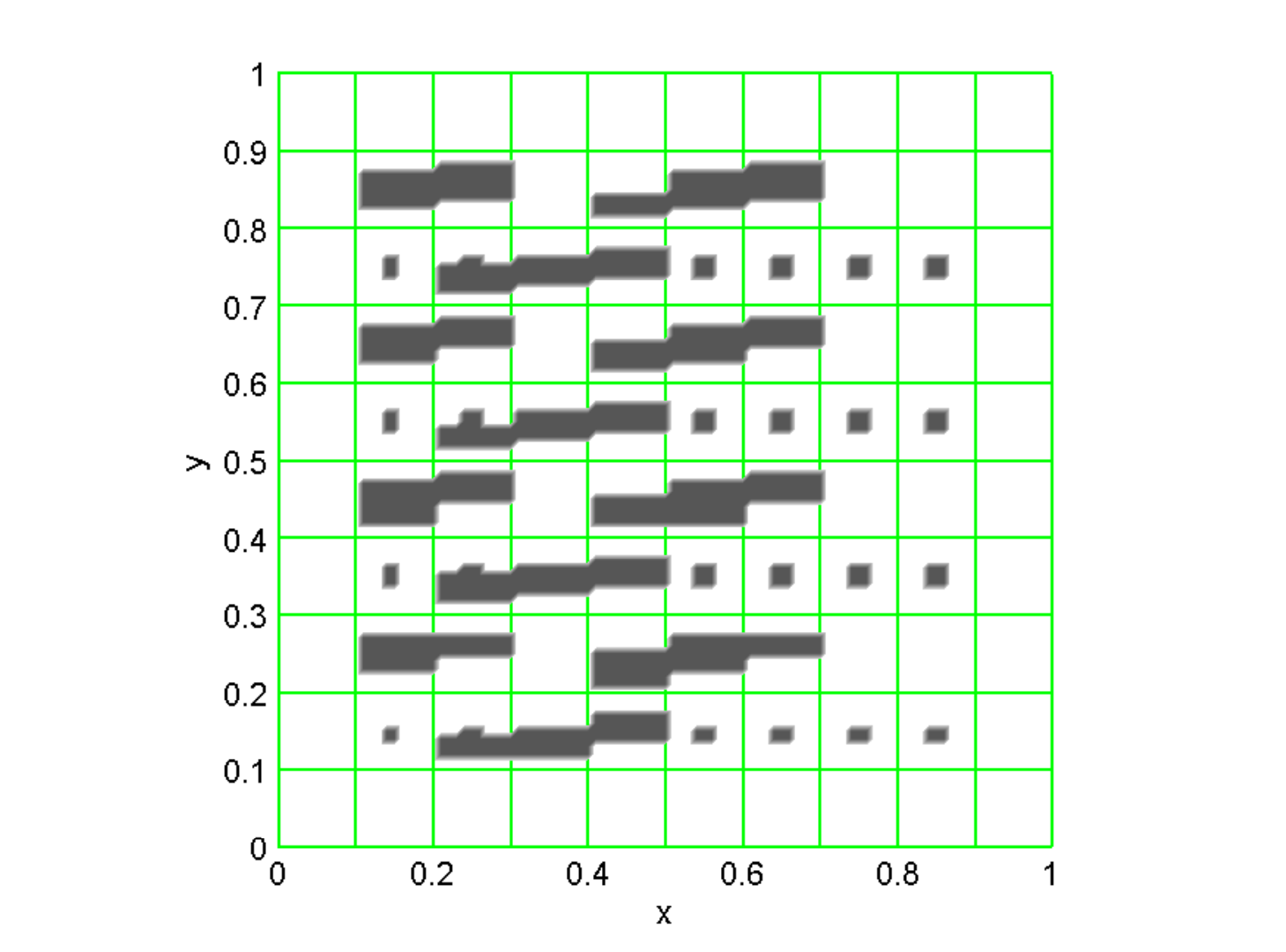}\label{Perm1}}
    \subfigure[Case
    II]{\includegraphics[width=0.48\textwidth]{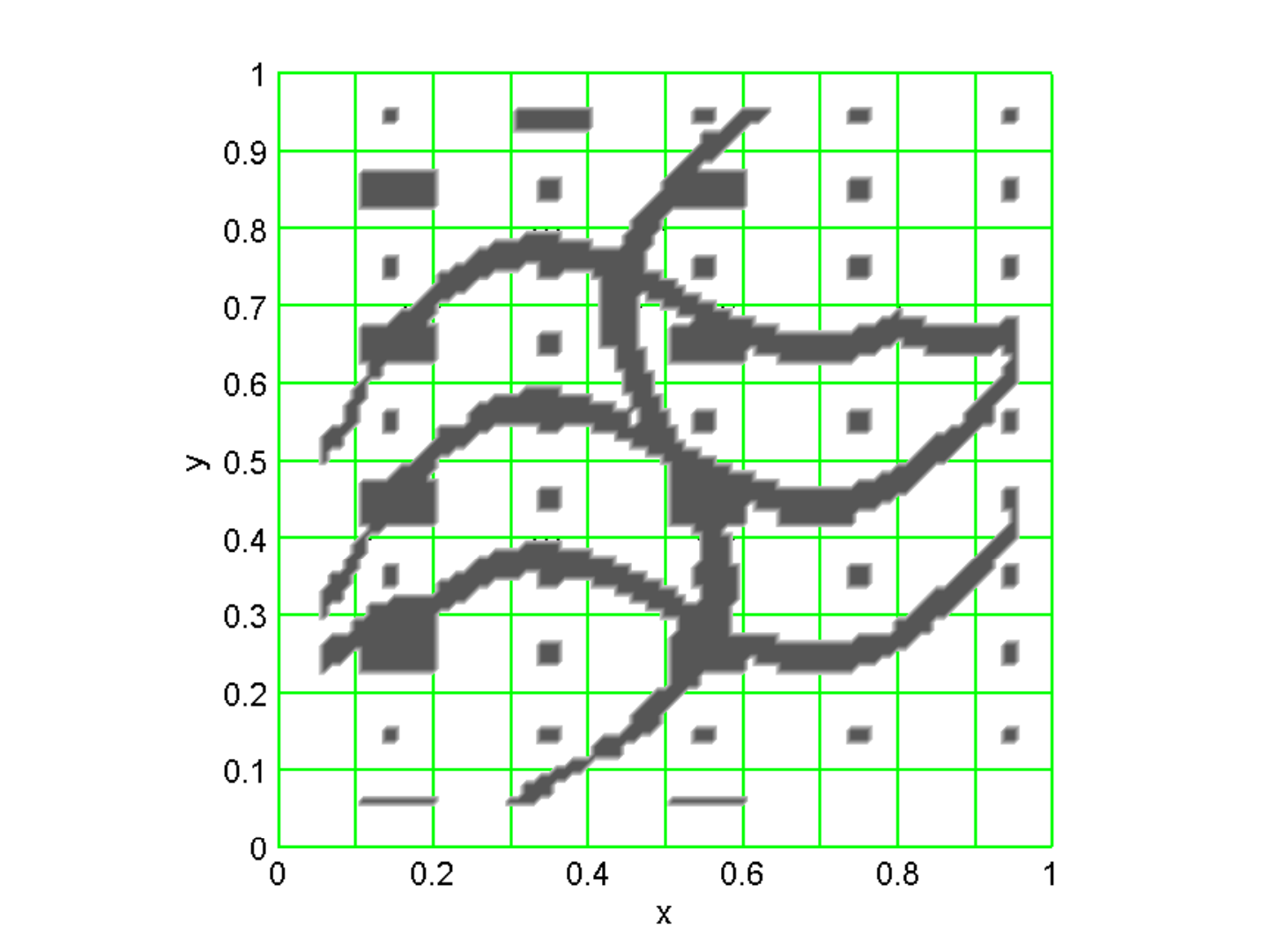}\label{Perm2}}
  \end{center}
  \caption{Different permeability fields considered in the numerical
    experiments.  Each permeability field has a different structures
    that model high conductivity channels within a homogeneous
    domain. The minimum (background) conductivity is taken to be
    $\kappa_{min} = 1$, and the high-conductivity (gray regions) with
    value of $\kappa_{max} = 10^{\eta} $ ($\eta$ = $4$ or $6$).}
  \label{Perm}
\end{figure}

In this section, we present numerical examples to illustrate the
applicability of the multiscale empirical interpolation method for
solving nonlinear multiscale partial differential equations. Before
presenting the individual examples, we review the computational domain
used in constructing the GMsFEM basis functions. This computation is
performed during the offline stage.  We discretize with finite
elements a nonlinear PDE posed on the computational domain
$D=[0,1]\times[0,1]$.  For constructing the coarse grid, we divide
$[0,1]\times[0,1]$ into $10\times 10$ squares. Each square is divided
further into $10\times10$ squares each of which is divided into two
triangles.  Thus, the discretization parameters are $1/100$ for the
fine-grid mesh and $1/10$ for the coarse-mesh. The fine-scale finite
element vectors introduced in this section are defined on this fine
grid. We recall that the fine-grid representation of a coarse-scale
vector $z$ is given by $\Phi z$, which is a fine-grid vector.

Using the two grids, we construct the GMsFEM coarse space as described
in Section~\ref{sec:app}. In the simulations, we consider two
different permeability fields.  These permeability fields, denoted by
$\kappa$, contain channels of high conductivity and are shown in
Figure~\ref{Perm}. This figure shows fields with different structures
that model high-conductivity channels within a homogeneous domain. The
minimum conductivity value for this case is taken to be $\kappa_{min}
= 1$, and the high conductivity varies from channel to channel with a
maximum value of $\kappa_{max} = 10^{\eta}$ ($\eta$ = 4 or 6).  Both
permeability fields are represented by the fine mesh.

\subsection{Nonlinear functions}

\begin{figure}[htb]
  \begin{center}
      \subfigure[$F_1$]{\includegraphics[width=0.47\textwidth]{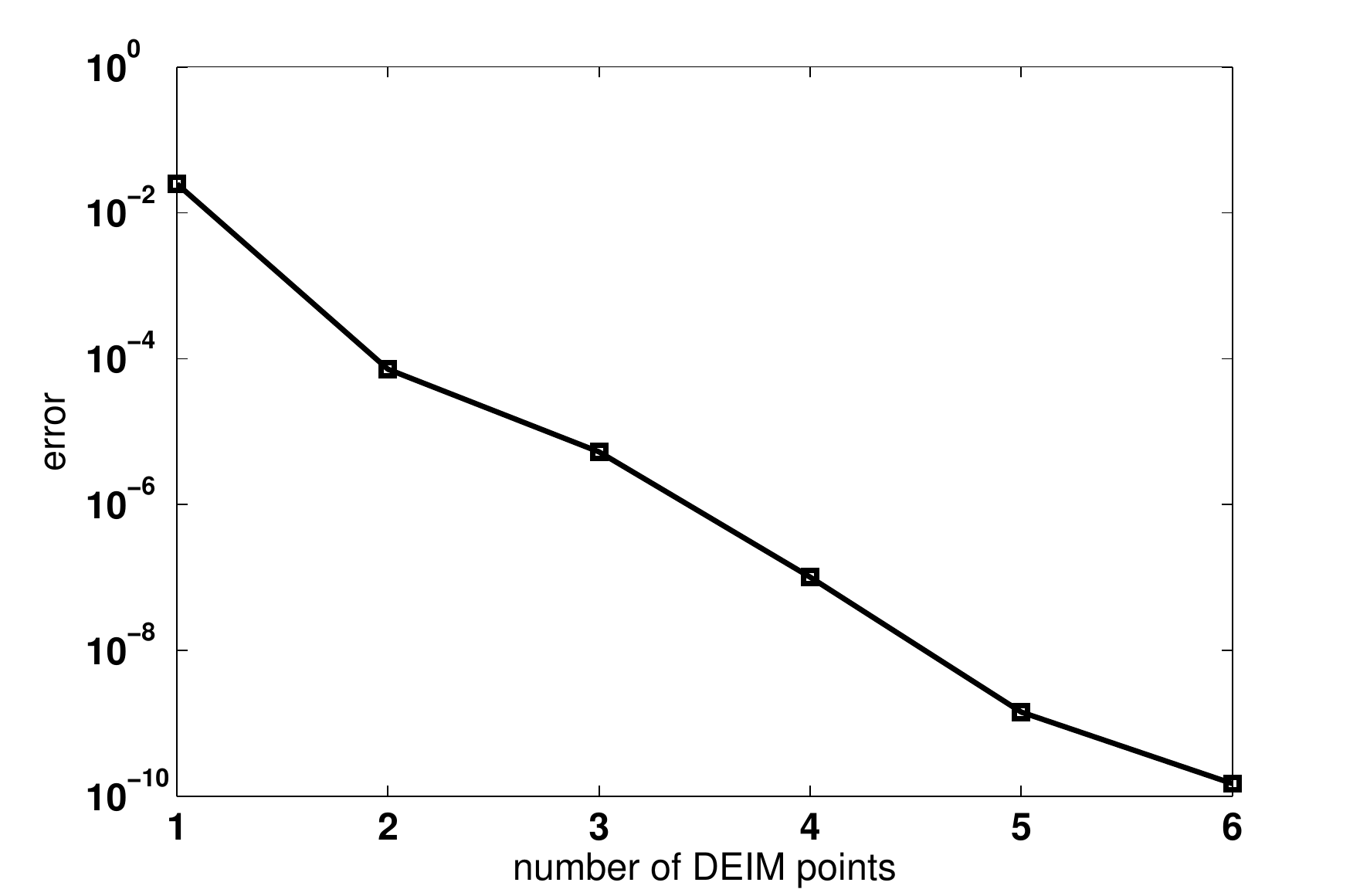}\label{errorF1}}
      \subfigure[$F_2$]{\includegraphics[width=0.47\textwidth]{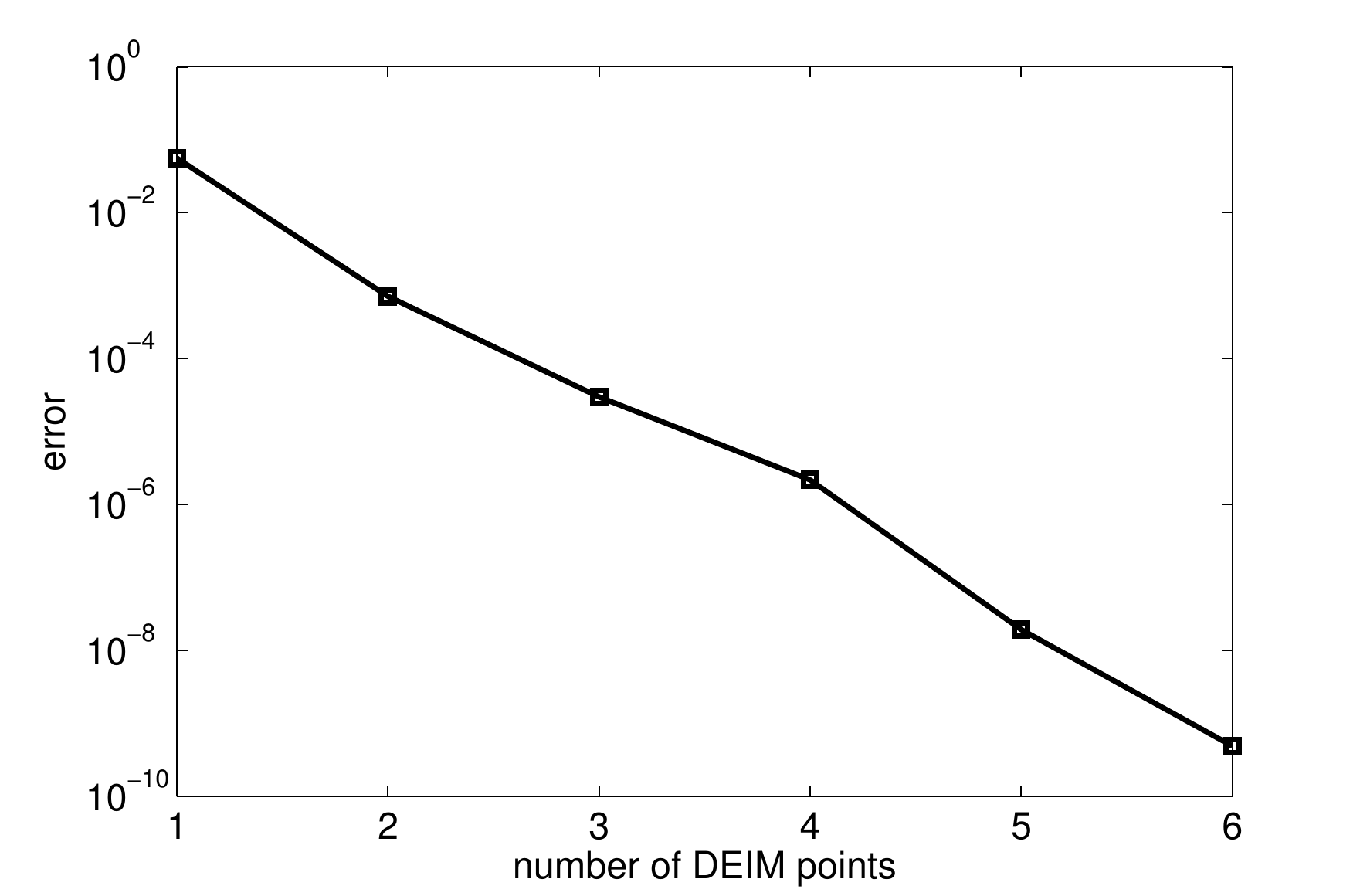}\label{errorF2}}
  \end{center}
  \caption{Variations of the relative error between the approximate
    and the original nonlinear functions with the number of DEIM
    points.}
  \label{Error1}
\end{figure}

To demonstrate the applicability of the multiscale empirical
interpolation method in approximating nonlinear functions when using
the multiscale framework, we consider two parametrized functions
$F_1:\mathbb{R}^{n_f}\times [0,1] \mapsto \mathbb{R}^{n_f}$ and
$F_2:\mathbb{R}^{n_f}\times [0,1] \mapsto \mathbb{R}^{n_f}$ given by
\begin{align}\label{NFs}
F_1(\Phi z;\mu) &= \Big(\sin(2\pi \mu \Phi z) \cos(2\pi \mu \Phi
z)\Big)^2e^{-2 \pi \mu \Phi z} \nonumber \\
F_2(\Phi z;\mu) &= \frac{1}{1+\sin(2\pi \mu \Phi z)},
\end{align}
where $\mu \in [0,1]$. We follow the multiscale DEIM approach
described in Section~\ref{MsDEIM} to approximate the nonlinear
function on the fine grid while evaluating at few selected points.

To check the capability of the multiscale DEIM to properly approximate
the nonlinear function, we compute the relative error as the
L$_2$-norm of the difference between the original and approximate
functions: i.e.,
\begin{eqnarray}\label{Err}
\|E\|_2 = \frac{\| F - \widetilde{F} \|_2}{\|\widetilde{F} \|_2},
\end{eqnarray}
where $\widetilde{F}$ is obtained from the DEIM approximation. We
consider $\Phi z$ as a multiscale solution of the elliptic
problem obtained from a fine mesh of dimension 10201
(i.e., $N_x=N_y=n_x=n_y=$10) by solving $-\mbox{div}(\kappa\nabla u)=1$
in $D$ and $u=0$ on $\partial D$, where
$\kappa$ is defined as in Case I with $\eta=4$ (see Figure \ref{Perm}).
 In Figure~\ref{Error1}, we plot the
relative error variations with the number of DEIM points.  As
expected, the error decreases as the number of DEIM points is
increased. For instance, the relative error $\|E\|_2$ is equal to 7.18
10$^{\mbox{-3}}$ and 7.14 10$^{\mbox{-3}}$ when using only 2 DEIM
points per region to approximate the nonlinear functions $F_1$ and
$F_2$, respectively, defined on a fine mesh of dimension 10201. These
results show the capability of the multiscale DEIM to reproduce the
fine-scale representation of the nonlinear function while gaining in
terms of computational cost through evaluating at only a few selected
points.

\begin{figure}[htb]
  \begin{center}
    \subfigure[$F_1$ - exact function
    ]{\includegraphics[width=0.45\textwidth]{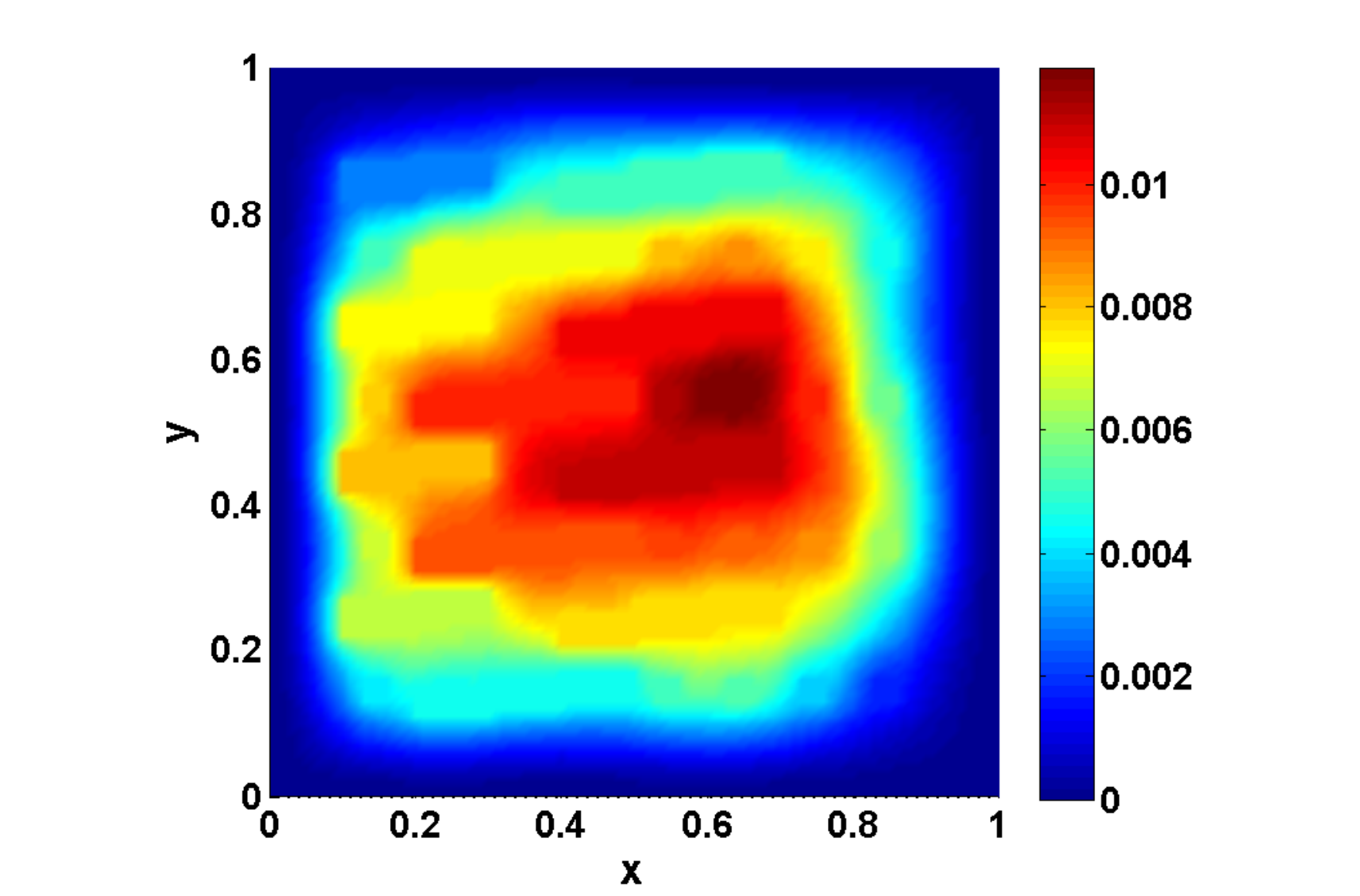}\label{F1}}
    \subfigure[$\widetilde{F}_1$ - approximate function
    ]{\includegraphics[width=0.45\textwidth]{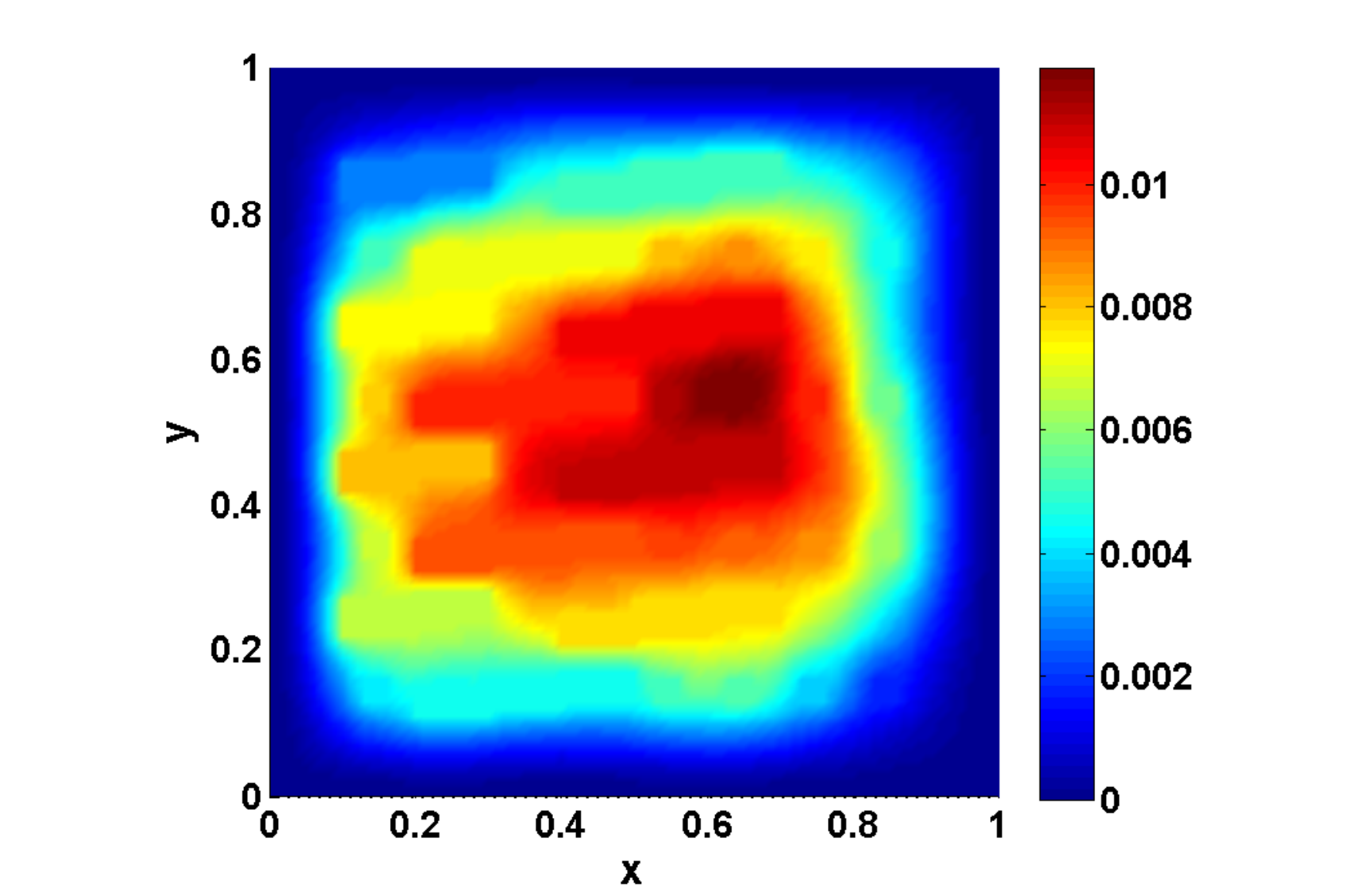}\label{F2}}
    \subfigure[$F_2$ - exact function
    ]{\includegraphics[width=0.45\textwidth]{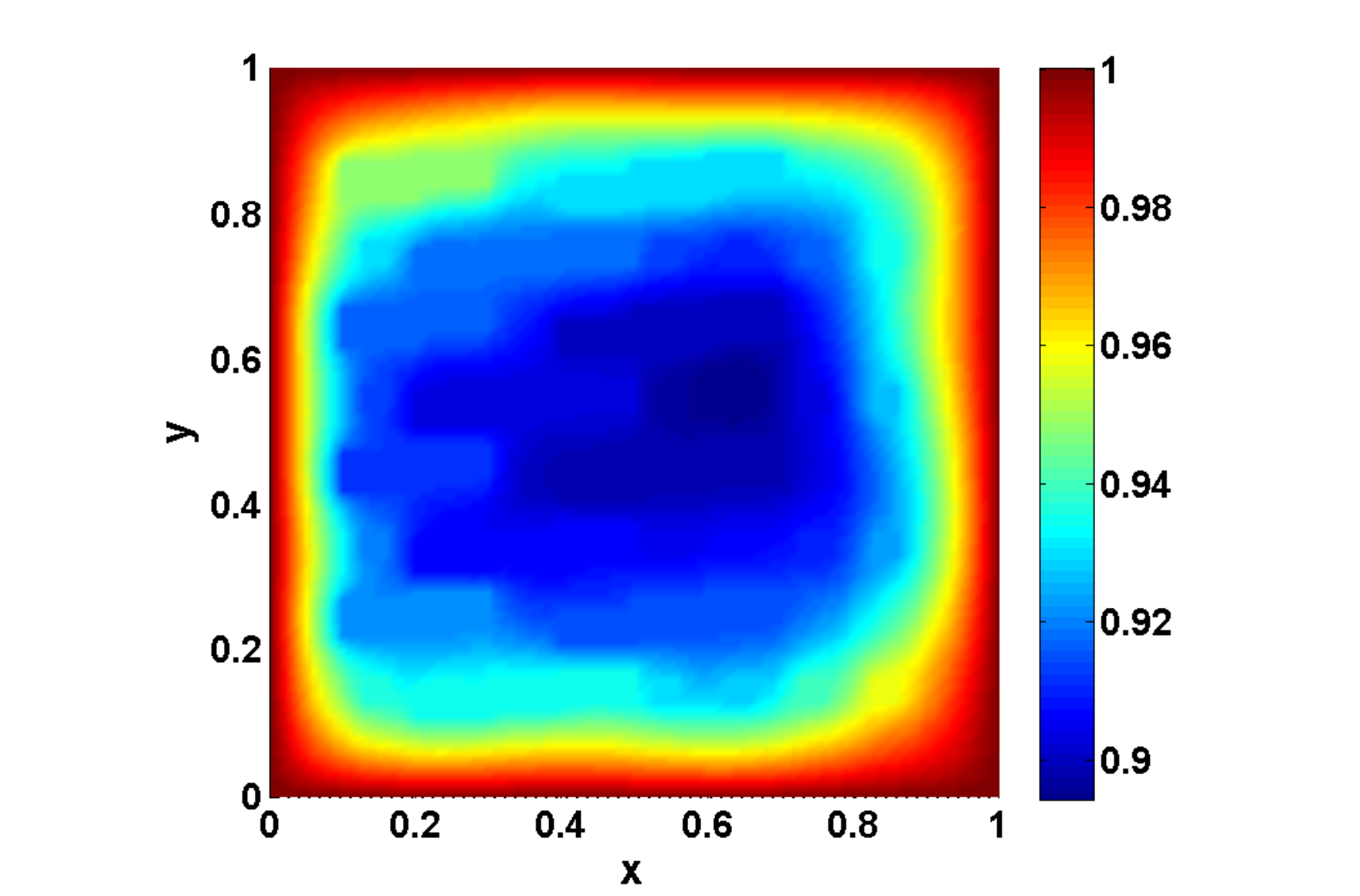}\label{F3}}
    \subfigure[$\widetilde{F}_2$ - approximate
    function]{\includegraphics[width=0.45\textwidth]{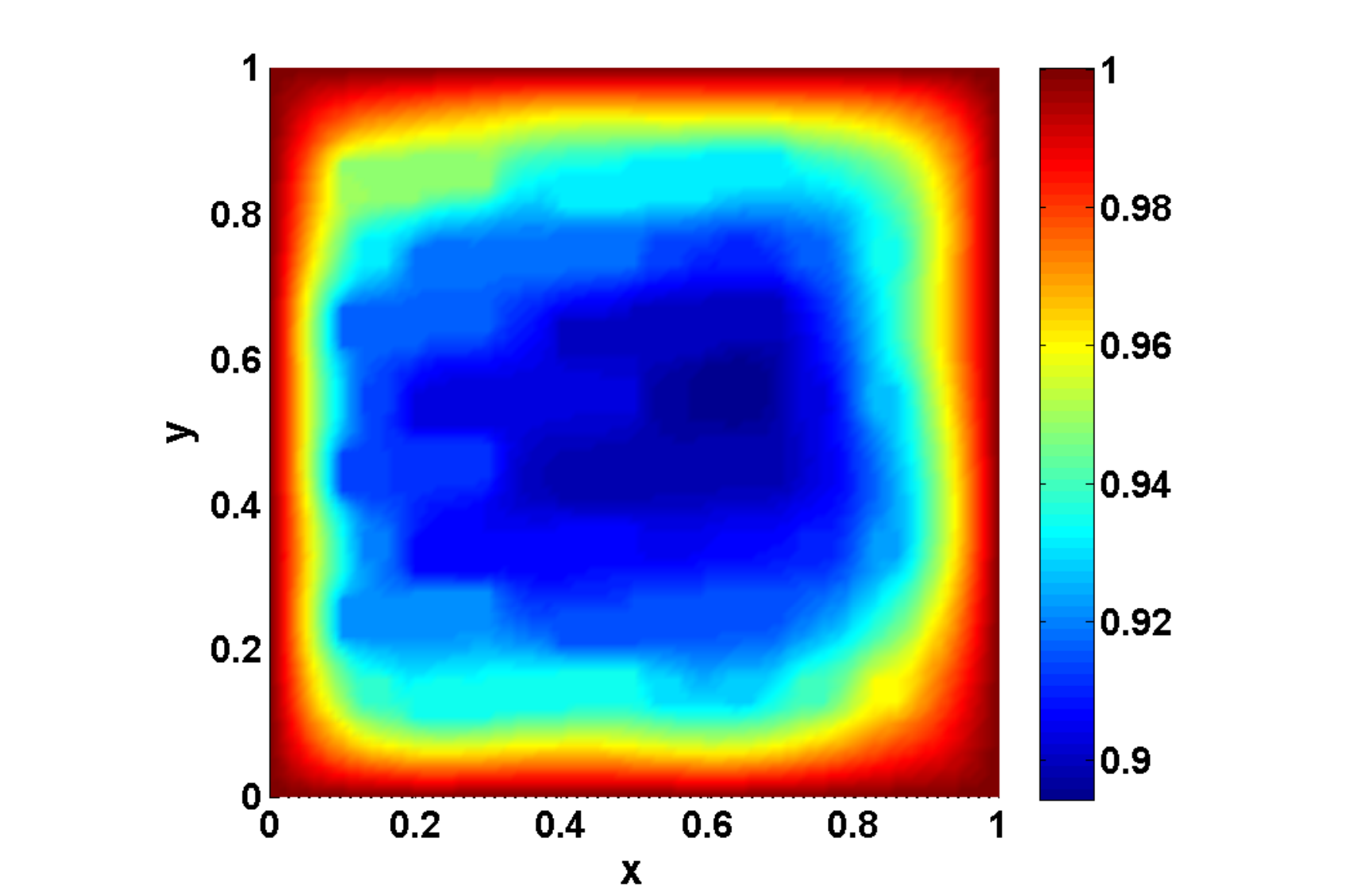}\label{F4}}
  \end{center}
  \caption{Comparison of the DEIM approximation with the original
    nonlinear function of dimension $n_f=10201$.}
  \label{F}
\end{figure}

In Figure~\ref{F}, we compare the approximate functions obtained from
DEIM using two points per region against the original function of
dimension $n_f=10201$. The good agreement observed between the two sets
of data shows the capability of DEIM to approximate the nonlinear
function using a few selected points per region.

\subsection{Nonlinear steady PDE}

\begin{figure}[htb]
  \begin{center}
    \subfigure[Case
    I]{\includegraphics[width=0.47\textwidth]{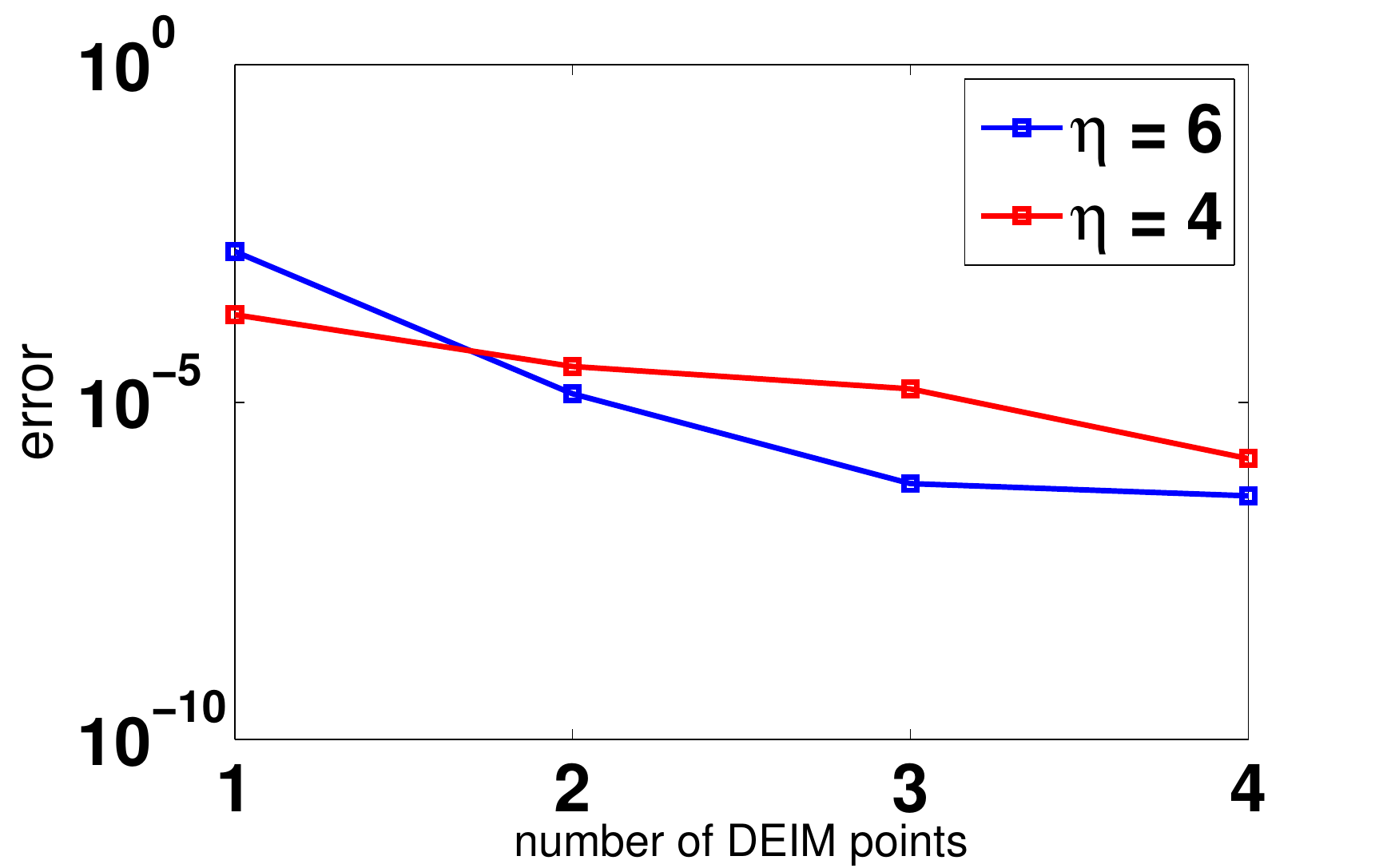}\label{errorPDE1Perm1}}
    \subfigure[Case
    II]{\includegraphics[width=0.47\textwidth]{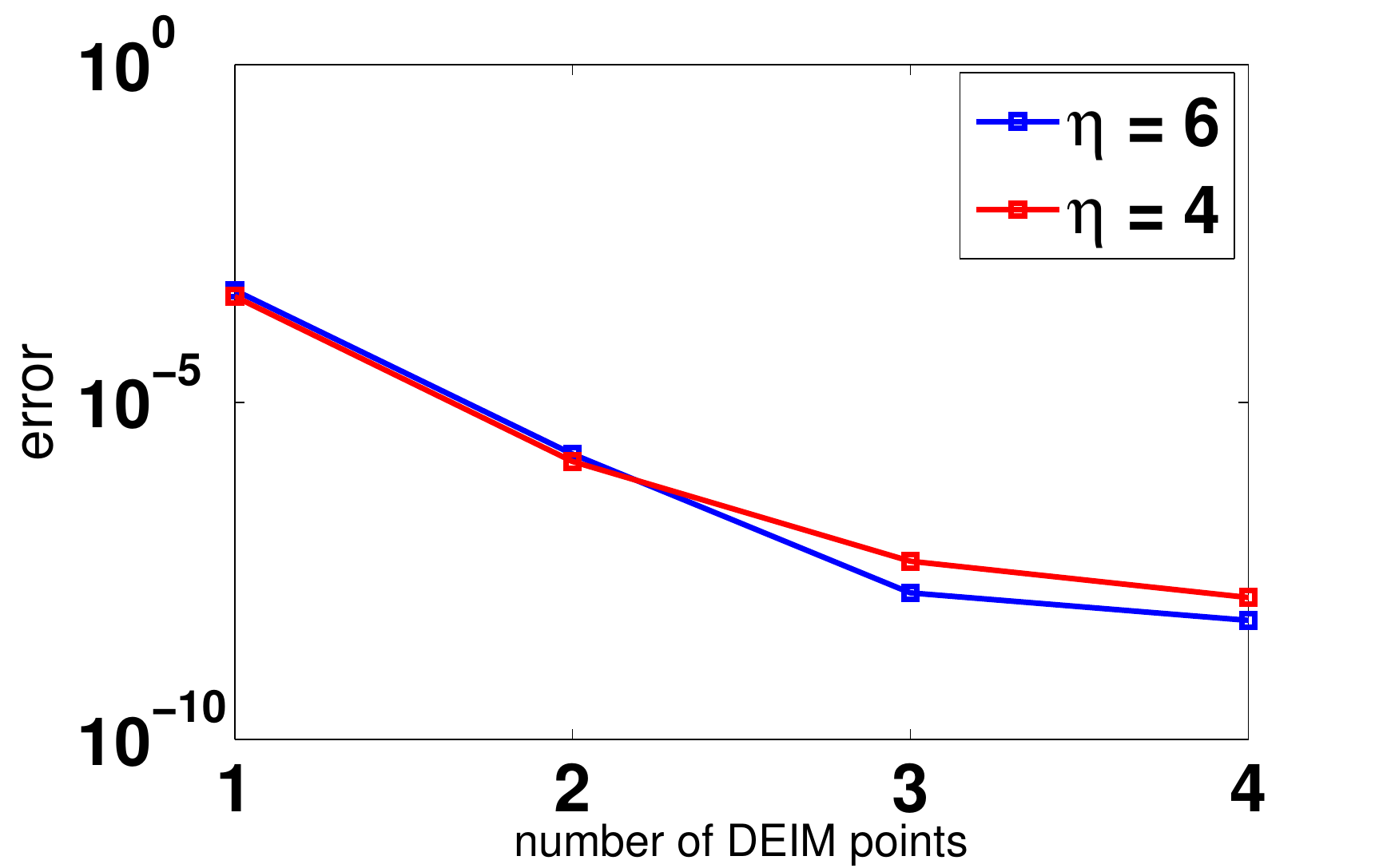}\label{errorPDE1Perm2}}
  \end{center}
  \caption{Variations of the relative energy error between the
    multiscale and multiscale DEIM solutions with the number of DEIM
    points.}
  \label{Error2}
\end{figure}

In this section, we consider a steady non-linear elliptic equation of the form
\begin{equation} \label{eqsteady}
\mbox{div} \big( \kappa(x) \, \nabla u  \big)=g(u,x,\mu) \, \,\, \, \text{in} \,\, \,\,  D,
\end{equation}
where
\begin{equation} \label{force}
g(u,x,\mu) = (1+\sin(2 \pi \mu u))e^{-2 \pi \mu u}.
\end{equation}
We consider permeability fields $\kappa$ that contain
high-conductivity channels as shown in Figure~\ref{Perm}.  We use
GMsFEM with Newton's method to discretize and solve
Equation~\eqref{eqsteady} and employ the multiscale DEIM to
approximate the nonlinear forcing term.

\begin{figure}[htb]
  \begin{center}
    \subfigure[Case
    I]{\includegraphics[width=0.47\textwidth]{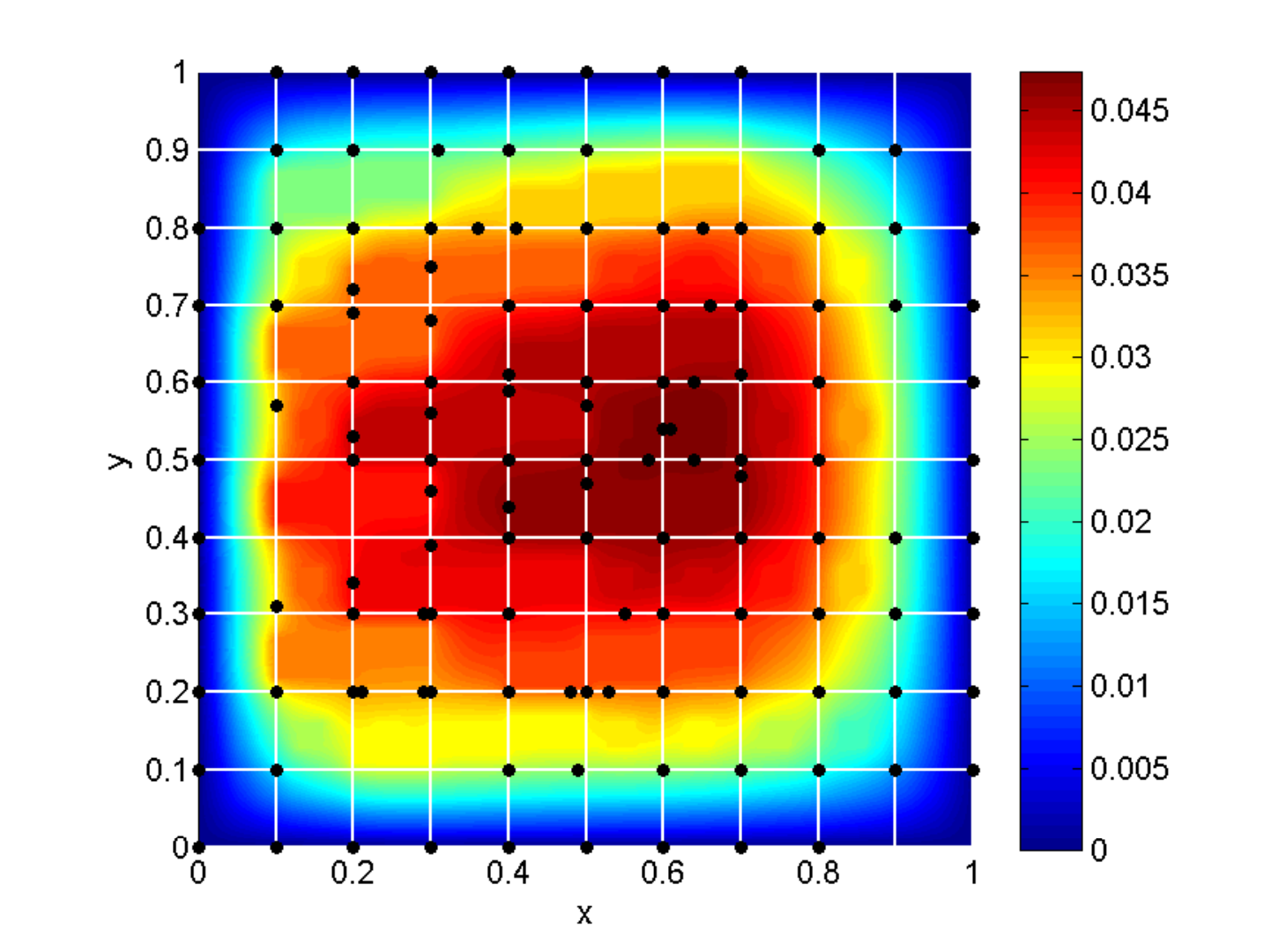}\label{DEIMPtsPDE1Perm1eta1}}
    \subfigure[Case
    II]{\includegraphics[width=0.47\textwidth]{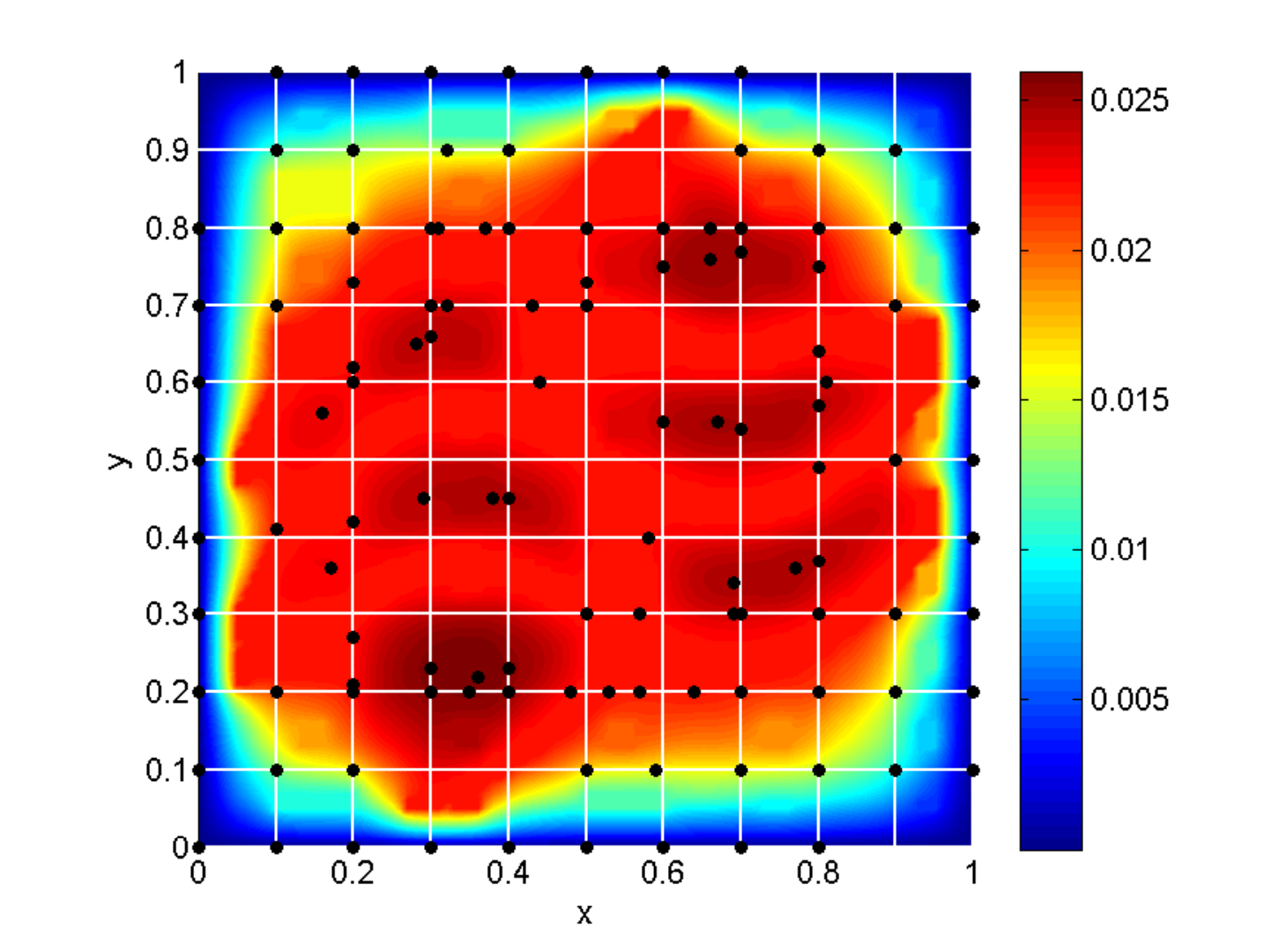}\label{DEIMPtsPDE1Perm2eta1}}
  \end{center}
  \caption{Coarse grid and approximated solutions obtained from the multiscale DEIM
    approach using two DEIM points per region. Black dots denote the
    location of the DEIM points.}
  \label{DEIMPDE1}
\end{figure}

We define the relative energy error as
\begin{eqnarray}\label{Err2}
  \|E\|_{\mbox{A}} = \sqrt{\frac{(\mbox{U} - \widetilde{\mbox{U}})^T
      {\mbox{A}} (\mbox{U} - \widetilde{\mbox{U}})}{\mbox{U}^T
      {\mbox{A}} \mbox{U}}},
\end{eqnarray}
where $\mbox{A}$ is the fine-scale stiffness matrix that corresponds
to~\eqref{eqsteady}. The errors are plotted in Figure~\ref{Error2}
with respect to the number of DEIM points. Clearly, the use of a few
DEIM points yields good approximation of the nonlinear forcing term
and, consequently, the nonlinear PDE solution. Our numerical
experiments show that the error does not depend on the contrast and it
decreases as we increase the number of DEIM points. We refer
to~\cite{g1,g2,egw10} for theoretical results on the error analysis of
high-contrast problems.

For illustration purposes, in Figure~\ref{DEIMPDE1}, we show the
approximate solutions obtained from the multiscale DEIM approach using
two DEIM points per coarse region. The black dots denote the location
of the DEIM points.  These DEIM points show high absolute values of
the modes representing the forcing terms in~\eqref{eqsteady}.

\subsection{Nonlinear unsteady PDE}

\begin{figure}[ht]
\begin{center}
\includegraphics[width=0.6\textwidth]{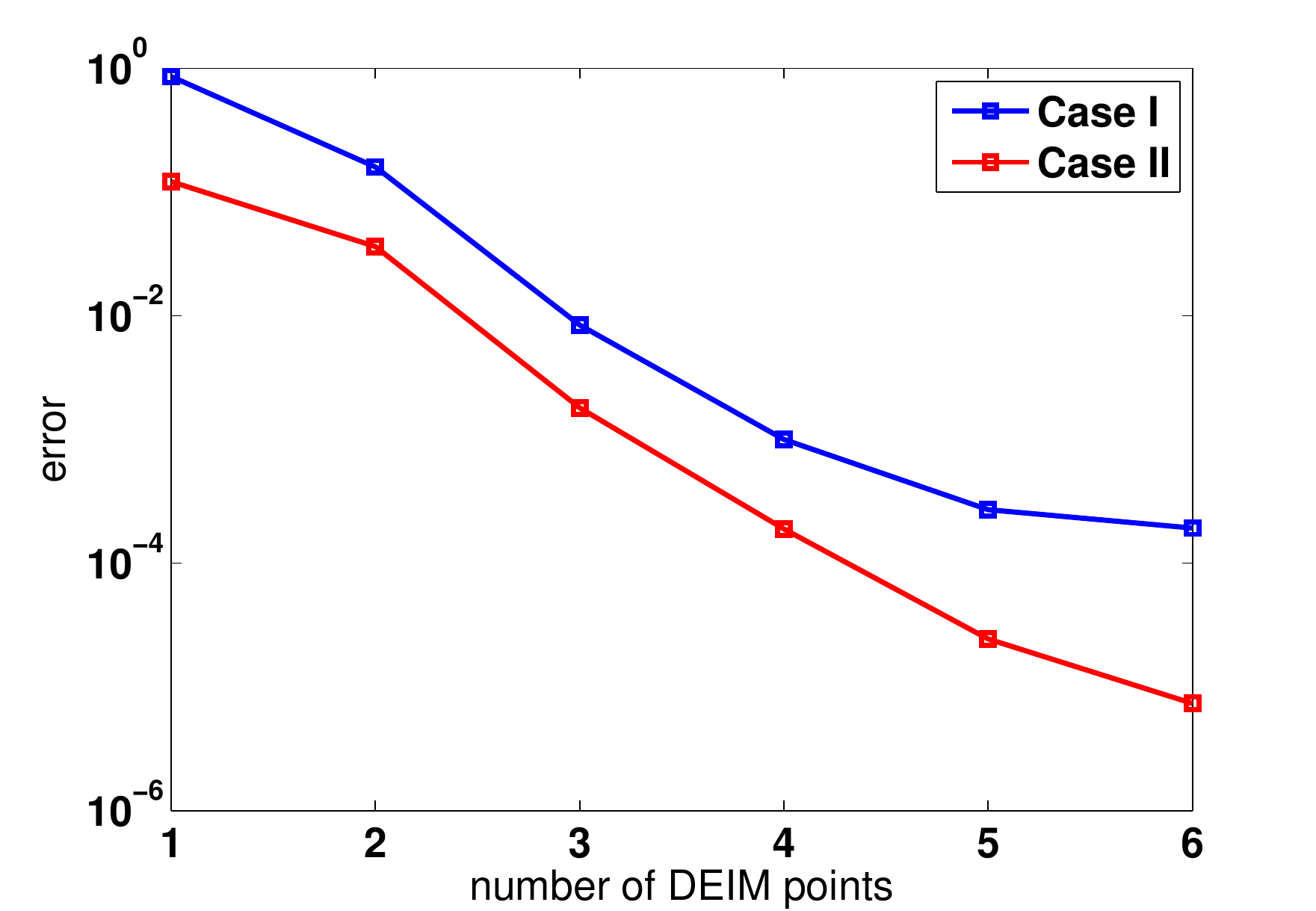}
\end{center}
\caption{Variations of the relative energy error between the
  multiscale and multiscale DEIM solutions with the number of DEIM
  points.}
  \label{Error3}
\end{figure}

As an example of a nonlinear unsteady problem, we consider the
following time-dependent parabolic equation
\begin{eqnarray}\label{parabolic2}
\frac{\partial u}{\partial t} -\mbox{div} \big( \kappa(x;u,\mu) \,
\nabla u  \big) = h(x)   \, \, \text{in} \, D,
\end{eqnarray}
where
\[
\kappa(x;u,\mu)= \kappa_q(x) b_q(u,\mu), \;\;\; b_q(u,\mu)=e^{\mu u},
\;\;\; h(x)=1+\sin(2\pi x_1)\sin(2 \pi x_2),
\]
and the structure of the permeability fields $\kappa_q$ is shown in
Figure~\ref{Perm}. We employ the GMsFEM for space discretization and
the Newton method to solve the nonlinear algebraic system at each time
step as detailed in Section~\ref{Newton}. Furthermore, we use
multiscale DEIM, described in Section~\ref{MsDEIM}, to approximate the
nonlinear functions and use multiscale POD (see
Section~\ref{sec:MsPOD}) to identify the modes in the empirical
interpolation.  Without multiscale POD, the Newton method may converge
slowly or may not converge when using $l_2$-inner-product-based
POD~\cite{volkwein05}. In Figure~\ref{Error3}, we plot the variation
of the relative energy error between the multiscale and multiscale
DEIM solutions with the number of DEIM points for the two
configurations of the permeability field. The results are obtained for
$\mu = 10$. Multiscale DEIM approximates well the multiscale solutions
as indicated by the small error values while significantly reducing
the number of function evaluations.

\begin{figure}[htb]
  \begin{center}
      \subfigure[Case I]{\includegraphics[width=0.47\textwidth]{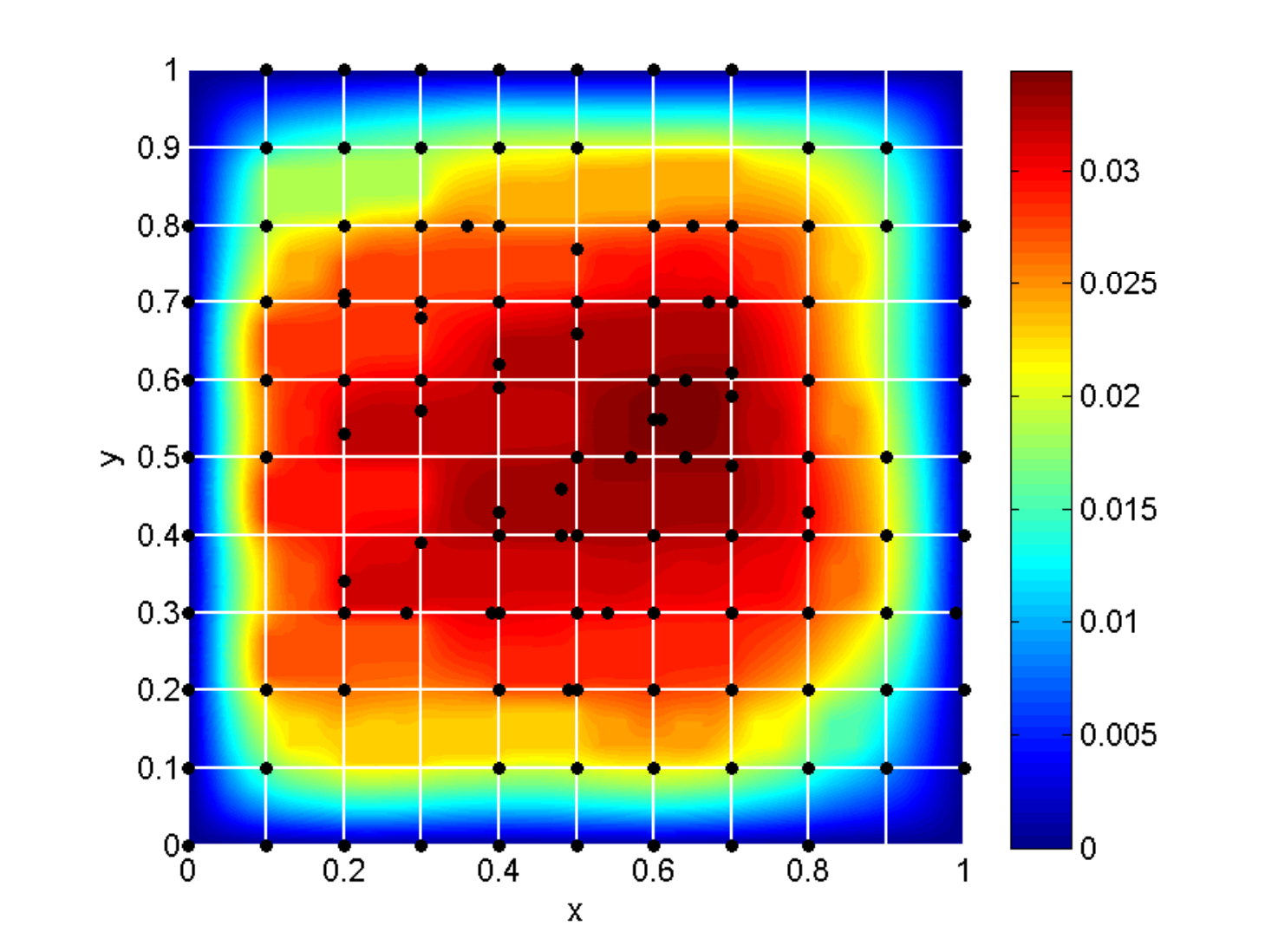}\label{DEIM2Perm1}}
      \subfigure[Case II]{\includegraphics[width=0.47\textwidth]{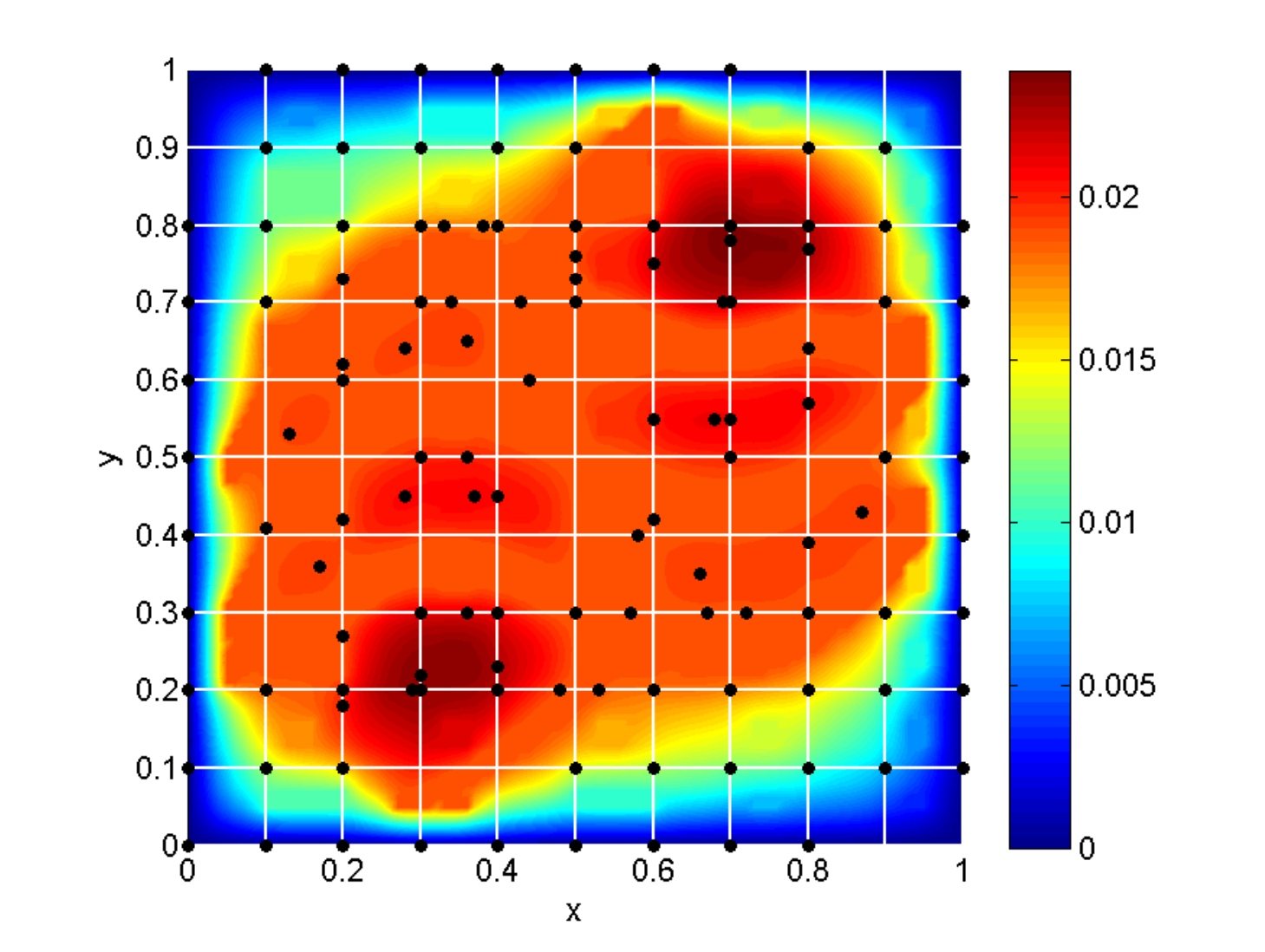}\label{DEIM2Perm2}}
  \end{center}
  \caption{Approximate solutions obtained from the multiscale DEIM approach using two DEIM points per region. Black dots denote the location of the DEIM points.}
  \label{DEIMPDE2}
\end{figure}

\begin{figure}[ht]
\begin{center}
\includegraphics[width=0.6\textwidth]{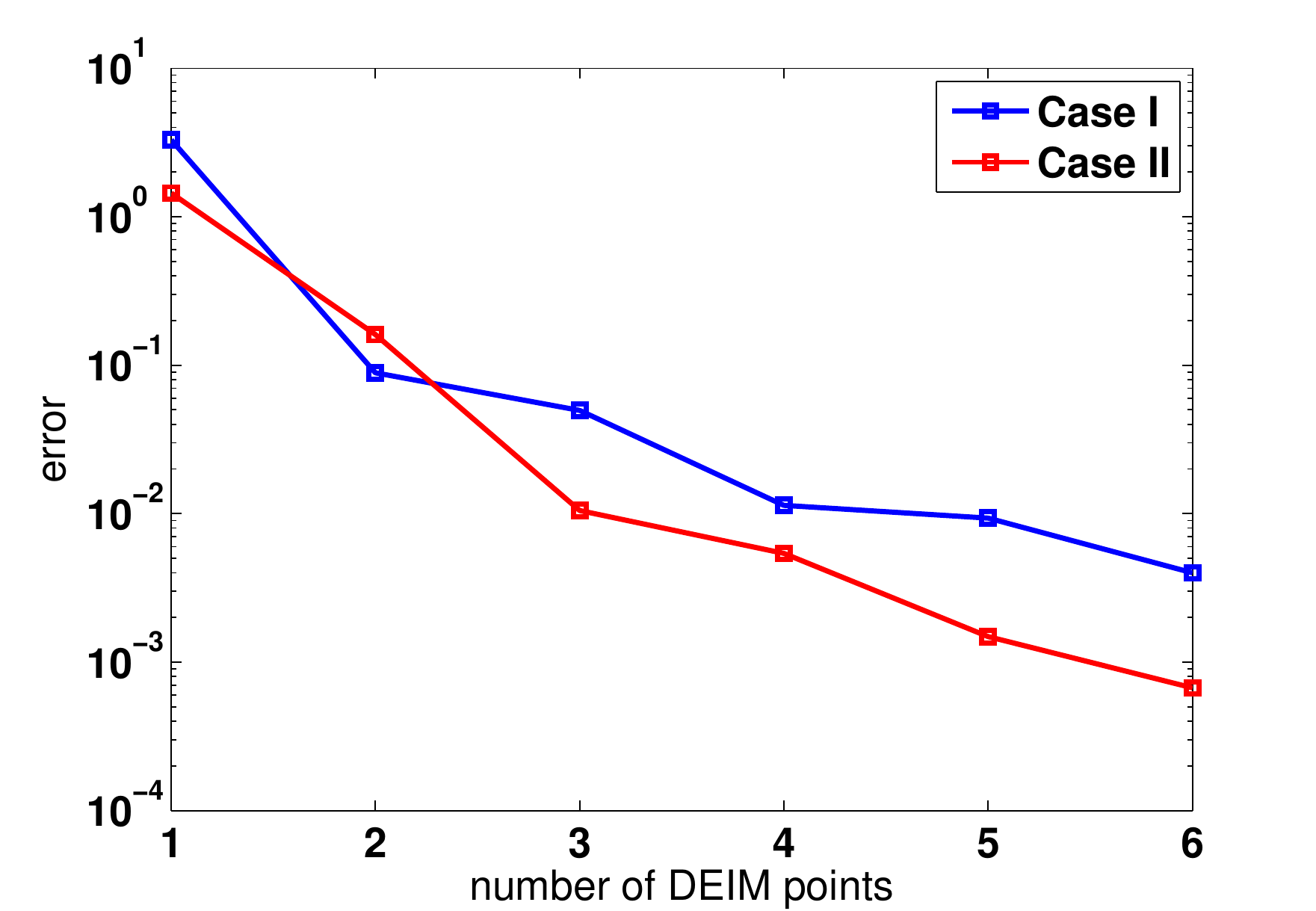}
\end{center}
\caption{Variations of the relative energy error between the
  multiscale and multiscale DEIM solutions with the number of DEIM
  points.}
  \label{Error4}
\end{figure}

In Figure~\ref{DEIMPDE2}, we show the approximate solutions obtained
from the multiscale DEIM approach using two DEIM points per
region. The black dots denote the location of the DEIM points.

We investigate the robustness of the multiscale DEIM approach with
respect to moderate variations in the nonlinearity of the permeability
field. As such, we run the forward problem while uniformly varying the
value of the parameter $\mu$ of the nonlinear function $b_q$ over the
interval $[5,10]$ with a step equal to one and store the snapshots of the nonlinear
functions for each case. Then, we use all snapshots to compute the
DEIM points and the index matrices of each region as described in
Section~\ref{MsDEIM}. Next, we consider a nonlinear function
$b_q(u,\mu)=e^{\mu u}$ with $\mu=8.5$ and show in Figure~\ref{Error4}
the relative energy error variations with the number of DEIM
points. Note that $\mu=8.5$ is not among the samples considered to compute the global DEIM points. Using a few DEIM points yields a good approximation.  We
mention that the simulation time of the final reduced order model is
about $4\%$ of the computational time of running the fine-grid model.
For instance, the relative energy error $\|E\|_{\mbox{A}}$ is equal to
4.96 10$^{\mbox{-2}}$ and 1.05 10$^{\mbox{-2}}$ when using only three
DEIM points per region to approximate the nonlinear functions that
appear in the residual $\widetilde{\mbox{R}}$ and the Jacobian
$\widetilde{\mbox{J}}$ for the two configurations of the permeability
field. Considering more DEIM points results in smaller error
values. These results demonstrate the robustness of the model
reduction approach based on combining the multiscale framework and
DEIM. We have also tested this approach with different
right-hand-sides and observed similar results.

\section{Conclusions}

In this paper, we propose a multiscale empirical interpolation for
solving nonlinear multiscale partial differential equations.  The
proposed method uses the Generalized Multiscale Finite Element Method
(GMsFEM) which constructs multiscale basis functions on a coarse grid
to solve the nonlinear problem. To approximate inexpensively nonlinear
functions that arise in the residual and Jacobians, we design
multiscale empirical interpolation techniques that use empirical modes
constructed based on local approximations of the nonlinear functions
using weighted POD techniques.  The proposed multiscale empirical
interpolation techniques (1) divide the computation of the nonlinear
function into coarse regions, (2) evaluate the contributions of the
nonlinear functions in each coarse region taking advantage of a
reduced-order solution representation, and (3) introduce multiscale
Proper Orthogonal Decomposition techniques to find appropriate
interpolation vectors.  We demonstrate the applicability of the
proposed methods on several examples of nonlinear multiscale PDEs that
are solved with Newton methods.  Our numerical results show that the
proposed methods provide an accurate and robust framework for solving
nonlinear multiscale PDEs on coarse grids while providing significant
computational cost savings..

\end{document}